\documentclass[12pt,leqno]{article}
\usepackage{amsmath,amsfonts,amsthm,amscd,amssymb,mathrsfs,esint}
\usepackage{xcolor,bbm,easybmat,bm}
\usepackage[colorlinks,citecolor=red,pagebackref,hypertexnames=false]{hyperref}
\usepackage{graphicx}
\usepackage{esint} 
\usepackage{dsfont}

\numberwithin{equation}{section}

\newtheorem{lem}{Lemma}[section]
\newtheorem{pro}[lem]{Proposition}
\newtheorem{defi}[lem]{Definition}
\newtheorem{thm}[lem]{Theorem}
\newtheorem{cor}[lem]{Corollary}
\newtheorem{rem}[lem]{Remark}

\newcommand{\ms}{\medskip}

\newcommand{\R}{\mathbb{R}}

\newcommand{\bY}{\mathbb {Y}}
\newcommand{\bT}{\mathbb {T}}
\renewcommand{\H}{\mathcal H}

\renewcommand{\d}{\partial}
\newcommand{\dist}{\,\mathrm{dist}\,}

\newcommand{\sm}{\setminus}

\renewcommand{\i}{\subset}
\newcommand{\wt}{\widetilde}

\newcommand{\cL}{{\mathcal L}}

\newcommand{\cP}{{\mathcal P}}
\newcommand{\W}{{\bf W}}
\newcommand{\X}{{\bf X}}
\newcommand{\cF}{{\cal F}}
\newcommand{\1}{{\mathds 1}}

\newcommand{\ol}{\overline}
\newcommand{\Per}{\mathrm{Per}}
\newcommand{\nn}{\nonumber}
\newcommand{\wh}{\widehat}

\begin{document}

\title{An optimal partition problem for the localization of eigenfunctions}

\author{
Guy David\footnote{
G. David was partially supported by the European Community H2020 grant GHAIA 777822,
and the Simons Foundation grant 601941, GD.
} 
\, and Hassan Pourmohammad  }

\newcommand{\Addresses}{{
  \bigskip
  \footnotesize

 Guy David, \textsc{Universit\'e Paris-Saclay, CNRS, 
 Laboratoire de math\'ematiques d'Orsay, 91405 Orsay, France} 
 \par\nopagebreak
  \textit{E-mail address}: \texttt{guy.david@universite-paris-saclay.fr} \\
  
  Hassan Pourmohammad, \textsc{Department of mathematics, Tarbiat Modares University, Tehran, Iran} 
  \par\nopagebreak
  \textit{E-mail address}: \texttt{h.pourmohammad@modares.ac.ir}
}}

\maketitle

\abstract{
We study the minimizers of a functional on the set of partitions of a domain $\Omega \subset \R^n$
into $N$ subsets $W_j$ of locally finite perimeter in $\Omega$, whose main term 
is $\sum_{j=1}^N \int_{\Omega \cap \d W_j} a(x) \d\H^{n-1}(x)$. Here the positive bounded function
$a$ may for instance be related to the Landscape function of some Schr\"odinger operator. We prove the existence
of minimizers through the equivalence with a weak formulation, and the local Ahlfors regularity and uniform rectifiability of the boundaries 
$\Omega \cap \d W_j$.
}

\tableofcontents

\section{Introduction}\label{S1}

The main theme of this paper is the elementary existence and regularity 
theory for partitions of a domain $\Omega$ that minimize a functional whose main term
is the integral of a given function $a$ on the total boundary of the partition; see \eqref{1a3}.

Its initial motivation comes from the localization of eigenfunctions 
for Shr\"odinger operators $\cL = - \Delta + {\cal V}$ on a domain 
$\Omega \subset \R^n$, in relation to the so-called landscape function
introduced in \cite{FM}. The landscape function (some times also called torsion function) 
is the solution $w$ of $\cL w = 1$ on $\Omega$, with the Dirichlet condition $w = 0$ 
on $\R^n \sm \Omega$, and the main point of \cite{FM} is that $w$ often controls the eigenfunctions
of $\cL$, both in the sense that there is a pointwise inequality 
of the type $u(x) \leq C_\lambda w(x)$ when $\cL u = \lambda u$,
but also because in many practical cases $w$ is localized in certain regions of $\Omega$,
and the eigenfunctions $u$, with not too high eigenvalues, also turn out to be localized in the same regions.

It seems interesting to find algorithms that, given $w$, find a collection of regions where 
$w$ is localized and hence it is hoped that the eigenfunctions live too. Of course it is tempting do this by 
minimizing a functional. Most often, in practical calculations, various forms of watershed 
algorithms have been used. The clear advantage of such methods is that they are often very fast, 
but one can argue that in some cases the description of valleys for $-w$ is perhaps not optimal.

The issue was taken in \cite{DFJM}, where a minimizing process was proposed involving the search for optimal free boundaries,
as in the Alt-Caffarelli-Friedman problem \cite{ACF}, and its regularity studied.  This gave
interesting practical results too, but the functional was not much used for the localization of eigenfunctions
since then, because the computations tend to be very long. 
Here we intend to use a different scheme that uses almost minimal partitions, i.e., relies on a functional whose main term measures a weighted version of the total area of the interfaces. 
This is simpler and more direct (and probably can be computed much faster) than the 
free boundary problem, but for some reason the authors of \cite{DFJM} were distracted and focused on 
free boundaries, and forgot about their initial attempt with surface measure.
The two approaches have common features, such as the use of many phases (many subdomains in
our partitions) and an auxiliary term (here called $G$).

\ms
The choice of $a$ and $G$ in relation with the localization of waves will be discussed more in Section \ref{S2},
but let us now describe the strong form of the (slightly more general) functionals that we intend 
to study in this paper.
We are given a domain $\Omega \i \R^n$, and an integer $N \geq 1$ (a bound for the number 
of pieces that we want). We want to use a functional $F$ (defined soon) to decompose 
$\Omega$ into disjoint subregions $W_{i}$, $1 \leq i \leq N$. Our set of acceptable competitors
will be the class $\cF = \cF(\Omega)$ of $N$-tuples 
$\W = (W_1, W_2, \ldots W_N)$, where the $W_i$, $1 \leq i \leq N$, are Borel subsets
of $\Omega$ such that
\begin{equation} \label{1a1}
\Omega = \bigcup_{i=1}^N W_i \ \text{ and } \ 
W_i \cap W_j = \emptyset \text{ for } i \neq j.
\end{equation}
That is, we consider partitions of $\Omega$  into $N$ Borel subsets, but we allow some of them to be empty. 
[We will systematically use bold letters to denote $N$-tuples.]
Then we want to minimize on $\cF$ the functional $J = J_s$ defined by 
\begin{equation} \label{1a2}
J_s(\W) = F_s(\W) + G(\W),
\end{equation}
where the main term $F_s$ involves the total weighted surface measure of the interfaces, 
and $G$ can be seen as a more stable bulk integral whose main goal may be to prevent the 
minimization by trivial solutions. For the definition of $F_s$, we
give ourselves a positive function $a$ on $\Omega$, and set
\begin{equation} \label{1a3}
F_s(\W) = \sum_{i=1}^N \int_{\Omega \cap \d W_i} a(x) \d\H^{n-1}(x),
\end{equation}
where $\H^{n-1}$ denotes the Hausdorff measure of co-dimension $1$ (think about surface measure)
and  $\d W_i$ is (usual) boundary of $W_i$. We could also set 
\begin{equation}\label{1a4}
\d(\W) = \Omega \cap \bigcup_{i=1}^N \d W_i
\end{equation}
 and use
\begin{equation} \label{1a5}
\wt F_s(\W) = \int_{\d(\W)} a(x) \d\H^{n-1}(x);
\end{equation}
as we shall see near Corollary \ref{c5a2}, this is roughly equivalent.

Concerning the bulk part $G(\W)$, we shall explain in Section \ref{S2} why some term like 
this is needed to avoid too trivial solutions, but the main properties that we need for our existence 
and regularity results are that that the variations of $G$ are smaller than the expected variations of $F_s$. 
We will settle on the following H\"older regularity in terms of the $W_j$ (in the $L^1$ norm): 
we will assume that there exists an exponent $\alpha \in (\frac{n-1}{n}, 1]$ and a constant 
$C \geq 0$ such that 
\begin{equation} \label{1a6}
\big|H(\W) - H(\W')\big| \leq C \sum_{i=1} \big|W_i \triangle W_i'\big|^\alpha,
\end{equation}
where $\W = (W_1, \ldots W_N)$ and $\W' = (W'_1, \ldots W'_N)$, we denote by $A \triangle B$ 
the symmetric difference $(A \sm B) \cup (B\sm A)$, and $|A|$ stands for the Lebesgue measure of 
the set $A$.

 The main point of \eqref{1a6} is even weaker than this: we want to make sure that when $\W$ 
 and $\W'$ are the same except on some ball of radius $r$, then
 $|H(\W) - H(\W')| \leq C r^{\beta}$ for some $\beta > n-1$, so that minimizers for $J_s$
 are what we shall call quasi- or almost minimizers for $F_s$ alone, which is enough for our regularity results.

The strong functionals $J_s$ and $\wt J_s$ are simple to define, but it is well known that 
existence results are then hard to establish directly, and that for this it is better to use weak functionals,
where the partitions are defined in terms of sets of finite perimeter (Caccioppoli sets), and the functional
uses BV (bounded variation) norms. So we shall introduce a related weak functional $J_w$,
prove existence results for $J_w$ (see Section \ref{S3}), and later on prove that this functional
is equivalent to our initial strong functionals, in the sense that modulo a little bit of cleaning for the minimizers 
of $J_w$, the weak and strong functionals have the same minimizers. We shall do this in Section \ref{S5},
after we prove the basic (but very useful) Ahlfors regularity results for the boundary associated to
minimizers of $J_w$ (see Section \ref{S4}).

After this, we will be allowed to use the equivalence of the weak and strong functionals 
and discuss their regularity equivalently in either setting. 

Returning to the existence result (Theorem \ref{t3a1} in Section \ref{S3}), we will get it
under the assumptions that $\Omega$ is bounded, $a$ is continuous and positive on $\Omega$,
and the bulk term $G$ is bounded and continuous for the $L^1$-distance of sets (see \eqref{3a4}).

Then the main estimate is probably the local Ahlfors regularity of the boundaries of minimizing partitions. 
Think about the strong functional to simplify the discussion, but the properties are proved first for the
minimizers of $J_w$. We will only prove local estimates in $\Omega$ (so as to get more flexibility if needed
and not worry about the regularity of $\Omega$ or uniform estimates on $a$), and assume the 
H\"older continuity of the bulk term $G$ (as in \eqref{1a6}; see the main assumption \eqref{3a4}),
and local lower and upper bounds on $a$; see \eqref{4a7}. Then we will get good local 
Ahlfors regularity bounds \eqref{4a8}  on the total (reduced) boundary (the union of the (reduced) 
boundaries of the $W_i$), that depend only on the parameters at hand, but not the number $N$ of pieces. 
See Theorem \ref{t4a1}.
We shall also prove Ahlfors regularity bounds on the individual boundaries of the $W_i$, but so far with bounds 
that depend on $N$; see Theorem~\ref{t4a4}.

The equivalence of weak and strong functionals will be proved in Section~\ref{S5}, after we 
describe how to clean a minimizer $\W \in \cF_w$ to replace the $W_i$ with equivalent open sets, whose boundaries
$\d W_i$ are contained, modulo an $\H^{n-1}$-negligible set, in the corresponding reduced boundaries $\d^\ast W_i$. 

Section \ref{S6} contains the next expected step, where we prove that the clean boundary $\d(\W)$ of a minimizer 
$\W$ is locally uniformly rectifiable, and even satisfies S. Semmes' Condition $B$; the same is even 
true for the individual boundaries $\d W_i^\sharp$ of the (equivalent open) pieces, but with an estimate that depends on $N$. See Theorem \ref{t6a1}.
With a little more work, we show that the $W_i^\sharp$ have a nice shape, in the sense that some local isoperimetric inequality is satisfied in $W_i^\sharp$. See Theorem \ref{t6a2}, and the discussion that follows on the fact that
this means a local John estimate (on the existence of thick curves that allow to escape from a point
of $W_i^\sharp$.

In Section \ref{S7} we check that if in addition $a$ is assumed to be (locally) H\"older continuous,
the clean boundary $\d(\W)$ of a minimizer is a local almost minimal set, as in \cite{Al} and \cite{Ta}.
This allows one to use the difficult regularity results obtained for such sets, and in particular, in dimension $n=3$, 
the beautiful local description of the $W_i^\sharp$ and $\d(\W)$ of \cite{Ta}, where $\d(\W)$ is shown to be
$C^{1+\varepsilon}$-equivalent to one of the three minimal cones that show up in soap films. 

We conclude the paper in Section \ref{S8} with a short list of open questions.

This paper is not the first one where functionals like $F_s$ and $F_w$
are studied; for instance, they were used in \cite{DS2} (but with only two pieces $W_i$), 
with similar uniform rectifiability 
results, to find big pieces of uniformly rectifiable sets inside sets that satisfy a topological condition
and a conflicting upper bound on Hausdorff measure. And in \cite{Le}, functionals like $F_w$ were studied, 
but more subtle because the price of the interface between $W_i$ and $W_j$ comes from functions
$a_{i,j}(x)$ that depend on the pair $(i,j)$. In this more complicated case the local regularity results on the individual 
$\d W_i$ are still unknown, even when $N=4$.

The authors want to thank Gian Paolo Leonardi for interesting discussions on 
\cite{Le} and the related problem of infiltration.

\section{More on the localization of waves}\label{S2}

In this section we justify the idea of using a functional like $J_s$ above to decompose a domain $\Omega$
into regions where the eigenfunctions of Shr\"odinger operators operator $\cL = - \Delta + {\cal V}$
will tend to be supported. We use the fact that these eigenfunctions are dominated by the landscape function
$w$, i.e., the solution $w$ of $\cL w = 1$ on $\Omega$, with the Dirichlet condition $w = 0$ 
on the boundary (see \cite{ADFJM1, ADFJM2}), and thus try to decompose $\Omega$ into nice regions 
$W_i$ where $w$ is large, or equivalently that look like valleys for the effective potential $\frac{1}{w}$. 
The paper \cite{DFJM} contains a description of one of the ways, apparently indirect, to find such a decomposition 
with the help of the free boundary problem. Here we propose a more classical and simpler approach, where we
try to obtain that the boundaries $\d W_i$ are both rather short and tend to pass through the places where
$w$ is small by forcing the decomposition to minimize a functional $J_s$ whose main term $F_s$
is of the form 
$F_s(\W) = \sum_{i=1}^N \int_{\Omega \cap\d W_i} a(x) d\H^{n-1}(x)$ (as in \eqref{1a3}),
where $a$ is a continuous functions that is small when $w$ is small. For instance, we could take
$a(x) = \delta + w(w)$, but we did not test options numerically and other choices could perform better.
The reason for adding a small constant $\delta > 0$ is that $w$ is positive on $\Omega$, but vanishes on
the boundary because of its definition with a Dirichlet boundary condition. We usually do not want
the region near $\d\Omega$ to play a special role in the computations, so it seems reasonable
not to look for trouble by allowing $a$ to be too small (we will see that our regularity estimates may deteriorate where $a$ is too small). On the other hand, the standard regularity results for elliptic operators say that
$w$ is H\"older-continuous, and requiring $a$ to be H\"older continuous will be enough for us, so we do not need to do special efforts here.

If we just minimize $F_s$, something stupid happens: the minimal partition is just the trivial one where
for instance $W_1 = \Omega$ and all the other $W_i$ are empty, so that all the interior boundaries
$\d W_i \cap \Omega$ are empty. This disturbed the authors of \cite{DFJM} for some time, until they realized
that it is easy to add another piece of functional $G(\W)$ that prevents this, by favouring partitions $\W$
for which the $W_i$ have roughly equal volumes. So for instance, we could take
\begin{equation} \label{2a1}
H(\W) = \sum_{j=1}^N h(|W_i|),
\end{equation}
where $h$ is a nice function that tends to be convex. This is what was added (for instance)
to the free boundary part of the functional in \cite{DFJM}. We could also use the 
function $w$ itself to replace the volumes $|W_i|$ in \eqref{2a1} with some weighted integrals like
$\int_{W_i} q(w(x)) dx$, for some carefully chosen function $q$. We shall not try to optimize such choices here; let us just observe that it is easy to ensure that they satisfy our mild estimates on $G$. 

Generally speaking, we did not do experiments and we shall not try to optimize the various choices in the 
definition of $J_s$, but let us observe that a nice possibility is to also allow a special
region $W_0$ in the partition, that we could call the black region, and where in fact $w$ is
quite small and we expect the eigenfunctions to be small too. This region could be treated differently
than the other ones, typically by replacing the contribution 
$\int_{\Omega \cap\d W_0} a(x) d\H^{n-1}(x)$
with a smaller one like $\delta \H^{n-1}(\Omega \cap \d W_0)$
(we still want the boundary to be
nice, so we put a little bit of surface measure in the functional anyway, but maybe this is not even needed because $\Omega \cap\d W_0$
 will contribute indirectly because it is covered by the other $\d W_j$).

We refer the reader to Section 2 of \cite{DFJM} for a more detailed description of our motivations in regard to the localization problem, and the role of $G(\W)$.

\section{The weak functional $J_w$}\label{S3}

We now introduce the weak form of our functionals $F_w$ and $J_w$.
We consider a fixed bounded domain $\Omega \subset \R^n$, and 
call $\cP(\Omega)$ the set of Borel subsets $W \subset \Omega$ that have a finite 
perimeter in $\Omega$. 
This means that the Borel function $\1_W$ (the characteristic function of $W$)
lies in $BV(\Omega)$, i.e., has bounded variation in $\Omega$. Or in other words, the vector-valued 
distribution $\nabla \1_W$ is a finite measure (in $\Omega$; we do not look at what happens on 
$\d \Omega$). We can write this measure as $\nabla \1_W = v d\mu_W$, where $v$ is a Borel vector
valued function and $\mu_W$ is a finite positive Borel measure on $\Omega$, often denoted by $|D \1_W|$ 
and called the total variation measure of the (distribution) derivative $D \1_W$. 

For any open set $U \subset \Omega$, we shall denote by
\begin{equation} \label{3a1}
\Per(W,U) = \mu_W(U) = \int_U |D \1_W|
\end{equation}
the perimeter of $W$ in $\Omega$. For this and other information on $BV$ functions and Caccioppoli sets,
we refer to \cite{AFP} or \cite{Gi}.

Since this will not cost us much, we shall also work with the set $\cP_{loc}(\Omega)$
of Borel sets $W \subset \Omega$ such that the restriction of $W$ to any open set $U \subset\subset \Omega$
(i.e., whose closure in $\R^n$ is a compact subset of $\Omega$) lies in $\cP(U)$. Those are also the sets $W$ 
such that $\nabla \1_W$ is a locally finite measure on $\Omega$; we can keep the same notation for 
$\mu_W$ and its description, even though $\mu_W$ is no longer a finite measure.

Let $N \geq 1$ be given. Denote by $\cF_w$ the class $N$-tuples
$\W = (W_1, \ldots, W_N)$, where each $W_i$, $1 \leq i \leq N$, is a
set of locally finite perimeter in $\Omega$, and in addition
\begin{equation} \label{3a2}
\1_W = \sum_{i=1}^N \1_{W_i}.
\end{equation}

We write the fact that the $W_i$ form a partition of $\Omega$ this way, to insist on the fact that now 
we shall identify two subsets of $\Omega$ that differ only by a set of measure $0$. 
That is, $\cF_w$ is in fact an equivalence class modulo negligible sets.
Our functionals $F_w$ and $J_w$ will be defined on $\cF_w$, which means in particular that 
$J_w(\W) = J_w(\W')$ when $\W = \W'$ modulo sets of Lebesgue measure $0$; this was not 
the case for $J_s$ and $F_s$. 

We are now ready to define our functionals $F_w$ and $J_w$ on $\cF_w$. 
Let $a : \Omega \to (0,+\infty)$ be a  positive  function on $\Omega$. 
For $\W \in \cF_w$, we systematically set $\W = (W_1, \ldots, W_N)$, and use the positive 
(locally finite)  measure $\mu_{W_i} = |D \1_{W_i}|$ as above. Then set
\begin{equation} \label{3a3}
F_w(\W) = \sum_{i=1}^N \int_{\Omega} a(x) d\mu_{W_i}(x)
\end{equation}
for $\W \in \cF_w$; notice that $F_w(\W)$ may be infinite, but it is well defined.

As for $G$, for this section we shall simply assume that it is defined on $\cF_w$ and continuous 
for the $L^1$ distance. That is, set
\begin{equation} \label{3a4}
\dist(\W,\W') = \sum_{i=1}^N \big|W_i \triangle W_i'\big|
\end{equation}
where $\W = (W_1, \ldots W_N)$, $\W' = (W'_1, \ldots W'_N)$, and we let
\begin{equation} \label{3a5}
W_i \triangle W'_i = (W_i \sm W'_i) \cup (W'_i\sm W_i);
\end{equation}
we only require that for $\W \in \cF_w$, 
\begin{equation} \label{3a6}
\lim_{k \to +\infty} G(\W_k) = G(\W)
\end{equation}
 when $\{ \W_k \}$ is a sequence in $\cF_w$ such that $\lim_{k \to +\infty} \dist(\W_k,\W) = 0$.
 Finally we set as expected
 \begin{equation} \label{3a7}
J_w(\W) = F_w(\W) + G(\W) \ \text{ for } \W \in \cF_w(\Omega).
\end{equation}

\begin{thm}\label{t3a1} 
Let $\Omega$, $a$, $G$, and $J_w$ be as above. Assume that $\Omega$ is bounded,
that $a$ is continuous and $a(x) > 0$ on $\Omega$, 
and that $G$ is bounded on $\cF_w$ and continuous for the distance of \eqref{3a4}.
Then we can find $\W_0 \in \cF_w(\Omega)$ such that
\begin{equation}\label{3a9}
J_w(\W_0) \leq J_w(\W)
\ \text{ for } \W \in \cF_w(\Omega).
\end{equation}
\end{thm}

This will be a standard application of the compactness properties of $BV$ functions.
We need the continuity of $a$ here (or at least its  lower semi-continuity) because we want
$F_w$ to be be lower semi-continuous.

Let $\Omega$, $a$, $G$, and $J_w$ be as in the statement, and set
$m = \inf_{ \W \in \cF_w} J_w(\W)$; notice that $m \leq J_w(\W_{00}) = G(\W_{00})< +\infty$, where 
$\W_{00} = (\Omega,\emptyset, \ldots, \emptyset)$ and because $F_w(\W_{00}) = 0$.
Let $\{ \W_k \}$ be a minimizing sequence in $\cF_w$.
We assumed that $G$ is bounded because we do not want any useless complication; 
 since $\Omega$ is bounded, this would easily follow if we required $G$ to be uniformly continuous 
for the distance of \eqref{3a4}, for instance. Because of this, we can find $M$ such that 
$-M \leq G(\W) \leq M$ for $\W \in \cF_w$, and now for $k$ large,
\begin{equation} \label{3a10}
F_w(\W_k) \leq J_w(\W_k) + M \leq m + 1 + M \leq J_w(\W_{00}) + 1 + M \leq 2M+1,
\end{equation}
where $\W_{00} = (\Omega, \emptyset, \ldots \emptyset)$ is the trivial competitor with no boundaries.
Thus we have a control on the boundaries. Write $\W_k = (W_{1,k}, \ldots, W_{N,k})$. 

For each ball $B \subset \Omega$ such that $2B \subset \Omega$, the continuity of $a$
gives a constant $\delta_B > 0$ such that 
\begin{equation} \label{3a8}
a(x) \geq \delta_B \text{ for } x\in B.
\end{equation}
Then for $k$ large,
\begin{multline} \label{3a11}
\sum_{i=1}^N \Per(W_{i,k}; B) 
\leq \big(\inf_{B} a\big)^{-1} \sum_{i=1}^N \int_{B} a(x) d\mu_{W_{i,k}}(x) 
\cr \leq \delta_B^{-1} F_w(\W_k) \leq \delta_B^{-1} (2M+1).
\end{multline}
This means the functions $\1_{W_{i,k}}$ are a bounded family in $BV(B)$.
These functions are also all bounded by $1$ pointwise, and by a standard consequence of the 
Poincar\'e inequality, they form a relatively compact set in $L^1(B)$ (for the norm). 
Because of this, we can extract a subsequence, that we still denote by
$\{ \W_k \}$, such that for every ball $B$ as above, functions $\1_{W_{i,k}}$ converge to a limit $g_i $
in $L^1(B)$. In addition, we can even find a subsequence that converges
pointwise, the limit $g_i$ does not depend on $B$, it is the characteristic function of some $W_i \subset \Omega$, 
and also $\sum_i \1_{W_i} = \1_\Omega$ almost everywhere. 

Next we use the lower semi-continuity of the $BV$ norm, which says that for each ball $B$ as above
and each $i$, $\Per(W_i;B) \leq \liminf_{k \to +\infty} \Per(W_{i,k};B)$. 
In particular, $W_i \in \cP_{loc}(\Omega)$ and hence $\W \in \cF_w$.
But also, since $a$ is continuous (lower semi-continuous would have been enough here), we also get that 
\begin{equation} \label{3a12}
\int_{B} a(x) d\mu_{W_i}(x) \leq \liminf_{k \to +\infty} \int_{B} a(x) d\mu_{W_{i,k}}(x).
\end{equation}
This is also true, for the same reason, when we replace $B$ by any open set which is compactly contained
in $\Omega$, and, by taking an increasing sequence of open sets that exhausts $\Omega$, with
$B$ replaced by $\Omega$. 

We still need to prove that $\W$ is a minimizer, i.e., that $J_w(\W) \leq m$. But we already know from 
\eqref{3a12} that $F_w(\W) \leq \liminf_{k \to +\infty} F_w(\W_k)$, and since $\Omega$ is bounded,
the fact that $W_{i,k}$ converges in $L^1(K)$ to $W_i$ for each compact set $K$ and each index $i$ implies
that $\dist(\W_k,W)$ tends to $0$ and, by assumption, 
$G(\W) = \lim_{k \to +\infty} G(\W_k)$. We sum the two estimates and get that
$J_w(\W) \leq \liminf_{k \to +\infty} F_w(\W_k) \leq m$, as needed.
\qed

\section{The reduced boundary is Ahlfors regular}\label{S4}

In this section, we introduce the reduced boundaries associated to an acceptable partition 
$\W = (W_1, \ldots, W_N) \in \cF_w$, and prove Ahlfors regularity results for them
when $\W$ minimizes $J_w$. We start with the definition and basic properties
of the reduced boundaries.

Each $W_i$ is a set of locally finite perimeter in $\Omega$, 
so it has a reduced boundary, that we shall denote by $\d^\ast W_i$ or simply $\d^\ast_i$, 
and which is a rectifiable subset of $\Omega$ (we still don't look at what happens on $\d \Omega$). 
We refer to \cite{AFP, Gi} for general information on reduced boundaries, but let us recall 
the most important features. First of all, the total variation measure $\mu_{W_i} = |D\1_{W_i}|$ 
is given by
\begin{equation} \label{4a1}
\mu_{W_i} = \H^{n-1}_{\vert \d^\ast_i},
\end{equation}
where $\H^{n-1}_{\vert \d^\ast_i}$ denotes the restriction to $\d^\ast_i$ of the Hausdorff measure.
More precisely, we know that for $\H^{n-1}$-almost every $x \in \d^\ast_i$,
$\d^\ast_i$ has an approximate tangent plane $P_i(x)$ at $x$, and there is a (measurable)
way to choose a unit normal vector $n_i(x)$ (in fact, pointing inwards) so that 
\begin{equation} \label{4a2}
\nabla \1_{W_i} = n_i(x) \H^{n-1}_{\vert \d^\ast_i}  = n_i(x) \mu_{W_i},
\end{equation}
where the gradient is in fact a vector valued distribution; we happily mix gradients and differentials
here. 

By definition of $\d^\ast_i$, we have a nice asymptotic description of $W_i$, as being 
close to a half space bounded by $P(x)$, and more precisely
\begin{multline} \label{4a3}
\lim_{r \to 0} r^{-n} \big|\big\{ y \in W_i \cap B(x,r) \, ; \, \langle y-x, n(x) \rangle \geq 0 \big\} \big|
\\
= \lim_{r \to 0} r^{-n} \big|\big\{ y \in B(x,r) \sm W_i \, ; \, \langle y-x, n(x) \rangle \leq 0 \big\} \big|
= 0,
\end{multline}
so that in particular $W_i$ and its complement both have Lebesgue density $1/2$ at $x$.
This information will be useful when we merge two sets $W_i$ locally,
and need to compute the effect on the $\d^\ast_i$. 

The $\d^\ast_i$ are our analogues of the $\d W_i$ for the strong functional, and their union 
\begin{equation} \label{4a4}
\d^\ast = \d^\ast(\W) = \bigcup_{i\in I} \d^\ast_i
\end{equation}
will be the first important object of study. Before we state our first Ahlfors regularity result,
let us say a bit more on the structure of $\d^\ast$ and the $\d^\ast_i$.

\begin{lem}\label{l4a1} 
Let  $\W \in \cF_w$ be given.
For $\H^{n-1}$-almost every $x \in \d^\ast$, $x$ lies in exactly two of the sets $\d^\ast_i$, 
$1 \leq i \leq N$. Consequently, the functional $F_w$ of \eqref{3a3} is also given by
\begin{equation} \label{4a5}
F_w(\W) = 2\int_{\d^\ast} a(x) d\H^{n-1}(x).
\end{equation}
\end{lem}

\ms
Indeed, it is classical (see \cite{Gi}) that the reduced boundary $\d^\ast_i$ coincides,
modulo a set of vanishing $\H^{n-1}$-measure, with the measure-theoretic boundary of $W_i$ (in $\Omega$),
which can defined for instance as the set $\d_i^\sharp$ of points of $\Omega$ where both $W_i$ and its 
complement have a positive upper Lebesgue density; 
the main ingredient for this is the fact that if $x\in \d_i^\sharp$, the Poincar\'e inequality implies 
that the lower density of $\mu_{W_i}$ at $x$ is positive.
Now if $x \in \d^\ast_i$, then the definition says that the density of $W_i$ at $x$ is precisely $1/2$,
which means that we can find $j \neq i$ such that $x\in \d_j^\sharp$, hence, $\H^{n-1}$-almost
surely, $x \in \d_j^\ast$. The density at $x$ of both sets $W_i$ and $W_j$ is $1/2$
(by \eqref{4a3}), and so $x$ does not lie in any other $\d_k^\sharp$ or $\d_k^\ast$.
So almost every $x \in \d^\ast$ lies in exactly two sets $\d^\ast_i$, and now
\eqref{4a5} is the same as \eqref{3a3}.
\qed

We now state our main Ahlfors regularity estimate.

\begin{thm}\label{t4a1} 
Let $\Omega$, $a$, $G$, and $J_w$ be as above. Assume that $\Omega$ is bounded,
and that $G$ is H\"older continuous on $\cF_w$ for some exponent $\alpha \in (\frac{n-1}{n},1]$
(and for the distance of \eqref{3a4}), which means that there exists $C_\alpha \geq 0$ such that 
\begin{equation} \label{4a6} 
\big|H(\W) - H(\W')\big| \leq C_\alpha \dist(\W, \W')^\alpha
\ \text{ for } \W, \W' \in \cF_w,
\end{equation}
as in \eqref{1a6}. Also let $\Omega' \subset \Omega$ be open, and assume that
we can find $\delta > 0$ such that 
\begin{equation} \label{4a7}
\delta \leq a(x) \leq 1 \ \text{ for } x \in \Omega'.
\end{equation}
Finally suppose that $\W \in \cF_w(\Omega)$ minimizes $J_w$ in the class $\cF_w(\Omega)$. 
Then there is a constant $C_{ar} \geq 1$, which depends only on $n$, $\delta$, $\alpha$, and $C_\alpha$, 
such that for $x \in \Omega' \cap \d^\ast$ and $0 < \min(1,\dist(x,\d \Omega'))$,
\begin{equation}\label{4a8}
C_{ar}^{-1} r^{n-1} \leq \H^{n-1}(\d^\ast \cap B(x,r)) \leq C_{ar} r^{n-1}.
\end{equation}
\end{thm}

\ms
We are mostly interested in $\Omega' = \Omega$, but we want to allow the case when
$a$ only satisfies \eqref{4a7} locally. The assumption that $a(x) \leq 1$ is essentially a renormalization
(we could also let $a$ tend to $+\infty$ in some places near $\d\Omega$, but then we could apply the theorem to a slightly different functional). Similarly, we take $r \leq 1$ to be able to compare the two pieces of the functional that have different homogeneities and get a constant that does not depend on $r$ (or $\W$ as long as it is a minimizer). 
We do not need $a$ to be continuous here (but then the existence may fail).

We shall see later that \eqref{4a8} is important, first because it is the pillar  of the comparison between 
the strong and weak functional, but also because its proof is the main initial ingredient for our other regularity results. 
Later in this section, we will prove a local Ahlfors regularity estimate for the $\d^\ast_i$ individually,
(and even later look for more precise regularity results for the $\d_i^\ast$), but the constants will depend on $N$. 

The proof below follows the same path as classical results in similar settings; see for instance 
\cite{DS2, Le}. We will see that somewhat weaker minimality properties of $\W$
(of almost- or quasi-minimality) are enough for the conclusion, which is why we will try now 
to single out the main features of the proof. 

We intend to test the minimality of $\W$ by constructing competitors $\X \in \cF_w(\Omega)$, 
which coincide with $\W$ outside of a given ball $B$ such that $\ol B \subset \Omega'$,
and which will be obtained by pouring some sets $B \cap W_i$ into other sets $W_\ell$.
Set $I = \{ 1, 2, \ldots, N \}$; this will be more convenient because we often want to cut $I$ into subsets.
Then let $B = B(x,r) \subset \Omega'$ be given, and also choose a proper subset $I_0$ of $I$ and,
for each $i \in I_0$, a target index $\ell(i) \in I_1 : = I \sm I_0$. 
We want to pour each $W_i \cap B$, $i \in I_0$, into the target piece $W_{\ell(i)}$, 
and more precisely we define a new competitor
$\X = (X_1, \ldots, X_N)$ by
\begin{equation} \label{4a9}
X_i = W_i \sm B \ \text{ for } i \in I_0
\end{equation}
and 
\begin{equation} \label{4a10}
X_\ell = W_\ell \cup \Big(\bigcup_{i \in I_0 ; \ell = \ell(i)} X_i \cap B \Big) \ \text{ for } \ell \in I_1.
\end{equation}

Note that the $X_j$ also form an almost-everywhere partition of $\Omega$, like the $W_i$,
and $\X \in \cF_w$ because we only moved the sets $B \cap W_i$ from one component to another one, 
and it is known that finite unions or intersections of sets of finite perimeters have the same properties.
We will need to check what happens to the various pieces of the reduced boundaries $\d^\ast_i$. 

We start with a description of the new total reduced boundary $\d^\ast(\X)$ associated to $\X$
in $B$. We start with its restriction to (the interior of) $B$. 
We are mostly interested in the case when $x \in \d^\ast_i \cap \d^\ast_j$ disappears from $\d^\ast$, 
i.e., no longer lies in any $\d^\ast X_k$, and this happens (modulo a $\H^{n-1}$-negligible set of exceptions) 
when $i \in I_0$ and $j = \ell(i)$, or symmetrically when $j \in I_0$ and $i = \ell(j)$, 
or when both $i$ and $j$ lie in $I_0$, and $\ell(i) = \ell(j)$.
No point of $\d^\ast(\X) = \bigcup_{k\in I} \d^\ast X_k$ appears in $B$, 
so, modulo a set of vanishing $\H^{n-1}$-measure,
\begin{equation} \label{4a11}
\d^\ast(\X \cap B) \subset \d^\ast \cap B\sm \bigg( \Big(\bigcup_{i \in I_0} \d^\ast_i \cap \d^\ast_{\ell(i)} \Big) \cup 
\Big(\bigcup_{(i,j) \in I_0^2 ;  j \neq i \text{ and } \ell(i) = \ell(j)} \d^\ast_i \cap \d^\ast_j \Big)\bigg).
\end{equation}
The other inclusion also holds (modulo a negligible set), but we don't need it so we leave its proof.

On $\Omega \sm \ol B$, the two sets $\d^\ast = \d^\ast(\W)$ and $\d^\ast(\X)$ coincide
(because we did not change anything in this open set), so we are left with the set $\d^\ast(\X) \cap \d B$ to analyse. 
For this it will be easier to assume that the following additional condition on $B$ is satisfied: for each $i \in I$,
\begin{multline} \label{4a12}
\text{$\H^{n-1}$-almost every point } x\in W_i \cap \d B(x,r)
\\
\text{ is a Lebesgue density point for $W_i$ (in $\Omega$).}
\end{multline}
For this property, we assume that we have chosen a representative $W_i$ for each $W_i$, so that 
$\sum_i \1_{W_i} = \1_\Omega$; then \eqref{4a12} makes sense. The set of radii $r>0$ for which  
\eqref{4a12} holds may depend on our representative, but anyway, for each choice of representatives $W_i$
and each choice of $x$, the set of $r \in (0, \dist(x,\d\Omega))$ for which \eqref{4a12} fails is negligible, 
by Fubini and because almost every point of $W_i$ is a Lebesgue density point for $W_i$.

Return to $\d^\ast$ and $\d^\ast(\X)$. Let $r$ be such that \eqref{4a12} holds. 
By assumption, $\H^{n-1}(\d^\ast \cap \d B) = 0$ (because the points of $\d^\ast_i$ are 
points of density $1/2$ for $W_i$, not Lebesgue points),
but of course our cutting procedure may introduce new boundaries for the $X_i$ on $\d B$. 
In fact, for $x\in \d B$ such that $x\in W_i$ and $x$ is a density point of $W_i$, this happens if and only if
$i \in I_0$, and then $x \in \d^\ast(X_i) \cap \d^\ast(X_{\ell(i)})$. 
That is, modulo a $\H^{n-1}$-negligible set,
\begin{equation} \label{4a13}
\d^\ast \cap \d B = \emptyset
\ \text{ and } \d^\ast(\X \cap \d B) = \bigcup_{i \in I_0} W_i \cap \d B,
\end{equation}
and more precisely 
\begin{equation} \label{4a14}
\d^\ast(X_i)\cap \d B = W_i \cap \d B \ \text{ for $i\in I_0$}
\end{equation}
and
\begin{equation} \label{4a15}
\d^\ast(X_\ell)\cap \d B = \bigcup_{i \in I_0; \ell(i) = \ell} W_i \cap \d B
\ \text{ for $\ell \in I_1$.}
\end{equation}

\ms
We shall return to this when we compare $F_w(\W)$ and $F_w(\X)$, but let us 
first deal with the $G$-terms and compare $G(\W)$ with $G(\X)$. 
Observe that since $\W$ and $\X$ coincide on the complement of $B$, \eqref{3a4} yields
\begin{equation} \label{4a16}
\dist(\W, \X) = \sum_{i=1}^N \big|W_i \triangle X_i\big|
\leq \sum_{i=1}^N  \big(|W_i \cap B| +  |X_i \cap B| \big) \leq 2|B| \leq C r^n,
\end{equation}
and by \eqref{4a3}
\begin{equation} \label{4a17}
|G(\X) - G(\W)| \leq C_\alpha \dist(\W, \X)^\alpha \leq C^\alpha C_\alpha r^{\alpha n}.
\end{equation}
Since $\W$ is assumed to be a minimizer,
\begin{multline} \label{4a18}
F_w(\W) =  J_w(\W) - G(\W) \leq J_w(\X) - G(\W) 
\\
= F_w(\X) + G(\X) - G(\W) 
\leq F_w(\X) + C^\alpha C_\alpha r^{\alpha n}.
\end{multline}
Let us interpret this as a local almost minimality property.

\begin{defi}\label{d4a3} 
Let $\gamma > 0$ and the constant $C_\gamma \geq 0$ be given. 
We say that $\W \in \cF_w(\Omega)$ is a local almost minimizer for $J_w$, associated to the 
gauge function $h(r) = C_\gamma r^{n-1}$, when for each ball $B = B(x,r) \subset \Omega$,
with $r \leq 1$, and each $\X \in \cF_w(\Omega)$ that coincides with $\X$ outside of $B$, 
\begin{equation} \label{4a19}
F_w(\W) \leq F_w(\X) + C_\gamma r^\gamma r^{n-1}.
\end{equation}
\end{defi}

We picked $\alpha > \frac{n-1}{n}$ precisely to make sure that $\gamma = \alpha n - (n-1) > 0$, 
so we just establishes in \eqref{4a18} that $\W$ is a local almost minimizer for $J_w$, 
associated to the gauge function $h(r) = C_\gamma r^{n-1}$, with $C_\gamma = C^\alpha C_\alpha$.
This is the only minimizing property that we need for Theorem \ref{t4a1} (and we shall no longer need $G$).

Next we compare $F_w(\W)$ and $F_w(\X)$, and for this the only feature of $F_w$ that we need
is the fact that $\delta \leq a(x) \leq 1$ in $\Omega'$, which contains $\ol B$ and hence
the symmetric difference between $\d^\ast(\W)$ and $\d^\ast(\X)$. That is,
\begin{eqnarray} \label{4a20}
F_w(\X)-F_w(\W) 
&=& 2\int_{\d^\ast(\X) \sm \d^\ast(\W)} a(x) d\H^{n-1}(x)
- 2\int_{\d^\ast(\W)\sm \d^\ast(\X)} a(x) d\H^{n-1}(x)
\nonumber \\
&\leq& 2\H^{n-1}(\d^\ast(\X) \sm \d^\ast(\W)) - 2\delta \H^{n-1}(\d^\ast(\W)\sm \d^\ast(\X))
\end{eqnarray}
and by \eqref{4a19} (we keep the notation $\gamma = \alpha n - (n-1)$ and 
$C_\gamma = C^\alpha C_\alpha$), 
\begin{eqnarray} \label{4a21}
2\delta \H^{n-1}(\d^\ast(\W)\sm \d^\ast(\X)) \leq 2 \H^{n-1}(\d^\ast(\X) \sm \d^\ast(\W))
+ F_w(\W) -F_w(\X)
\nn \\
\leq 2 \H^{n-1}(\d^\ast(\X) \sm \d^\ast(\W)) + C_\gamma r^\gamma r^{n-1}
\end{eqnarray}
or, with $C'_\gamma = (2\delta)^{-1} C_\gamma$, 
\begin{equation} \label{4a22}
\H^{n-1}(\d^\ast(\W)\sm \d^\ast(\X)) \leq  \delta^{-1}\H^{n-1}(\d^\ast(\X) \sm \d^\ast(\W))
+ C'_\gamma r^\gamma r^{n-1}.
\end{equation}
This is a form of quasiminimality of $\W$ for the Hausdorff measure, and we shall not need more than 
this for the proof of  Theorem \ref{t4a1} (or Theorem \ref{t4a4} later).
What we will do is pick various choices of balls $B$, sets $I_0 \subset I$, and mappings 
$\ell : I_0 \to I_1 = I \sm I_0$, and compute the two sides of \eqref{4a22} to get valuable information.

\ms 
For the proof of the upper bound in \eqref{4a8}, we decide to pour all of $B$ into $W_1$, say.
That is, we take $I_1 = \{ 1 \}$, $I_0 = I \sm I_1$, and take $\ell(i) = 1$ for $i \in I_0$.
We get a competitor $X$ such that $\d^\ast(X) \cap B(x,r) = \emptyset$
(because in $B(x,r) \subset X_1$, so there is no boundary there), and of course $\d^\ast(X) \cap \d B(x,r)
\subset \d B(x,r)$ trivially. Then 
$\H^{n-1}(\d^\ast(\W)\cap B(x,r)) \leq \H^{n-1}(\d^\ast(\W)\sm \d^\ast(\X)) \leq C \delta^{-1} r^{n-1}$
by \eqref{4a22}, as needed.

For the lower bound in \eqref{4a8}, we start with a ball $B_0 = B(x,r_0) \subset \Omega'$,
we suppose that $r_0 \leq r_{00}$, where the small constant $r_{00}$ will be chosen soon, and that 
\begin{equation} \label{4a23}
\H^{n-1}(\d^\ast  \cap B_0) < \varepsilon r_0^{n-1}
\end{equation} 
for some small enough $\varepsilon$ (to be chosen soon), and we want to show by induction that
for $k \geq 1$, if we set $r_k = 2^{-k}r_0$ and $B_k = B(x,r_k)$, we also have 
\begin{equation} \label{4a24}
\H^{n-1}(\d^\ast  \cap B_k) < s^k \varepsilon r_k^n
\end{equation}
with we shall choose $s = 2^{-\gamma} \in (0,1)$. Then we will produce a contradiction.

We start with $r_0$, and choose $\rho \in (\frac{9r_0}{10}, r_0)$ 
such that the security condition \eqref{4a12} holds for $\rho$. This is easy to do
because \eqref{4a12} holds for almost every $\rho$.
Set $D = B(x, \rho)$, take a second copy $D'$ of $D$, and glue the two copies along
$\d D$ to get a topological sphere $\wt D = D \cup D'$. We do this so that we can apply Poincar\'e's
inequality to $\wh D$ without worrying about boundaries. 

Next pick a proper subset $I_0$ of $I$, and set $W = \bigcup_{i \in I_0} D \cap W_i$.
Notice that $W$ is a set of finite perimeter, in $D$ as well as in $\wh D$, and
the union $\wh W$ of $W$ and its copy $W' \subset D'$ is also of finite perimeter.
Notice that $\d^\ast W \cap D \subset \d^\ast$, we have a similar inclusion for
$\d^\ast W' \cap D'$, and on $\d D$ we have no contribution from
$\d^\ast \wh W$, because by \eqref{4a12} almost every point of $\d D$ is a point
of density $1$ or $0$ for $W$ (in $\Omega$), hence also for $\wh W$
(for $\wh D$, which we see locally as a tube around $\d D'$, maybe with a slightly different
metric). Thus $\wh W$ has no perimeter on $\d D$, and 
\begin{equation} \label{4a25}
\H^{n-1}(\d^\ast \wh D) = 2 \H^{n-1}(\d^\ast D) \leq 2 \H^{n-1}(\d^\ast \cap B_0)
\leq 2 \varepsilon r_0^{n-1}.
\end{equation}
The Poincar\'e inequality in $\wh D$ (proved exactly as in a ball or a sphere) yields
\begin{multline} \label{4a26}
\min(|W|, |D\sm W|) = \frac12 \min(|\wh W|, |\wh B\sm \wh W|) 
\\
\leq C \Per(\wh W; \wh B)^{\frac{n}{n-1}} = C \H^{n-1}(\d^\ast \wh D)^{\frac{n}{n-1}}
\leq C \varepsilon^{\frac{n}{n-1}} r_0^n
\end{multline}
by \eqref{4a25}. We choose $I_0$ so that $|W| \leq \frac12 |B|$ (otherwise, pick the complement
$I \sm I_0$), but $|W|$ is as close to $\frac12 |D|$ as possible. Then the minimum in \eqref{4a26}
is $|W|$, and we get that
\begin{equation} \label{4a27}
|W| \leq C\varepsilon^{\frac{n}{n-1}} r_0^n \leq C \varepsilon^{\frac{n}{n-1}} |D|.
\end{equation}
Recall that $|W| = \sum_{i \in I_0} |W_i \cap D|$; if $\varepsilon$ is small enough,
\eqref{4a27} says that there is no $W_i$ such that 
$\frac{1}{10} |D| \leq |W_i \cap D| \leq \frac{9}{10} |D|$, and hence (choosing $I_0$ as above)
all the numbers $|D|^{-1} |W_i \cap D|$ are very small, except exactly one, which is very close to $1$.
Let us assume that the large one corresponds to $i = 1$, and now \eqref{4a26} or \eqref{4a27} (with 
$I_0 = I \sm \{ 1\}$) means that
\begin{equation} \label{4a28}
|D \sm W_1| = \sum_{i > 1} |W_i \cap D| \leq C \varepsilon^{\frac{n}{n-1}} r_0^n.
\end{equation}

We claim that we can choose a new radius $r$, with $r_1 < r < \rho$, such that \eqref{4a12} 
holds for $r$ and 
\begin{multline} \label{4a29}
\H^{n-1}(\d B(x,r) \sm W_1)  
\leq (\rho-r_1)^{-1} \int_{r_1 < r < \rho}
\int_{\d B(x,r) \sm W_1} d\H^{n-1} dr
\\
= C (\rho-r_1)^{-1} \big|(B(x,\rho) \sm B(x,r_1))\sm W_1 \Big|
\leq  C r_0^{-1} |D-W_1|  \leq C \varepsilon^{\frac{n}{n-1}} r_0^{n-1}
\end{multline}
where the first inequality is possible by Chebyshev, the second one is a 
baby version of the co-area formula where we replace a bulk integral by an integral on spheres,
and then we use \eqref{4a28}. The additional property \eqref{4a12} does not disturb, because it holds 
for almost all $r$.

We use $B = B(x,r)$ to apply the construction of a competitor $\X \in \cF_w$ described above, 
with $I_0 = I \sm \{ 1 \}$. That is, all the small components are poured into $W_1$. 
Notice that with this choice, $\d^\ast(\X) \cap B$ is empty (because now only $X_1$ meets $B$), 
and so
\begin{multline} \label{4a30}
\H^{n-1}(\d^\ast(\X) \sm \d^\ast(\W)) \leq \H^{n-1}(\d^\ast(\X) \cap \d B)
\\
\leq \H^{n-1}(\d B \sm W_1) \leq C (s^k \varepsilon)^{\frac{n}{n-1}} r_k^{n-1}
\end{multline}
by \eqref{4a14} and \eqref{4a15}. Then by \eqref{4a22}
\begin{multline} \label{4a31}
\H^{n-1}(\d^\ast \cap B(x,r_1)) 
\leq \H^{n-1}(\d^\ast \cap B) \leq \H^{n-1}(\d^\ast(\W) \sm \d^\ast(\X)) 
\\
\leq C \H^{n-1}(\d^\ast(\X) \sm \d^\ast(\W)) + C'_\gamma r_0^\gamma r_0^{n-1}
\leq C \varepsilon^{\frac{n}{n-1}} r_0^{n-1} + C'_\gamma r_0^\gamma r_0^{n-1}.
\end{multline}
That is, $r_1$ satisfies \eqref{4a23} with $\varepsilon$ replaced by
$\varepsilon_1 = C \varepsilon^{\frac{n}{n-1}} + C'_\gamma r_0^\gamma$.

We are now ready to prove \eqref{4a24} by induction. First choose $\varepsilon$ so small
that $C \varepsilon^{\frac{1}{n-1}} \leq \frac14$. When $n=1$, we don't need to do a choice; in fact
we could choose $r$ so that the two points of $\d B$ lie in $W_1$, and then the term 
$C \varepsilon^{\frac{n}{n-1}}$ does not even exist; hence the strange power.
Now choose $r_{00}$ so small that $C'_\gamma r_{00}^\gamma \leq \frac14 \varepsilon$.
With these choices done, assume that \eqref{4a24} holds for $k$ (this is the case when $k=0$)
and prove it for $k+1$. By induction assumption, and the same argument as above but starting
from $r_k$ and $B_k$, we get that $r_{k+1}$ satisfies \eqref{4a23}, with
\begin{multline}\label{4a32}
\varepsilon_{k+1} = C (s^k \varepsilon)^{\frac{n}{n-1}} + C'_\gamma r_k^\gamma
\leq C \varepsilon^{\frac{1}{n-1}} s^k \varepsilon  + C'_\gamma 2^{-k \gamma} r_0^\gamma
\\
\leq \frac14 s^k \varepsilon + \frac14 2^{-k \gamma} \varepsilon
\leq s^{k+1} \varepsilon 
\end{multline}
because we took $s = 2^{-\gamma} \geq 2^{-1}$. Thus \eqref{4a24} holds for all $k \geq 0$,
and in particular 
\begin{equation}\label{4a33}
\lim_{\rho \to +\infty} \rho^{1-n}\H^{n-1}(\d^\ast \cap B(x,\rho)) = 0.
\end{equation}
We claim that this is impossible for $x \in \d^\ast$. Indeed for $x \in \d^\ast_i$, 
\eqref{4a3} says that for $\rho$ small enough, we can find two balls $B_1, B_2 \subset B(x,\rho)$,
of radius $\rho/4$, such that if $m_{i,j}$ denotes the average of $\1_{W_i}$ on $B_j$, $j = 1, 2$,
$m_{i,1}$ is as close to $1$ as we want and $m_{i,2}$ is as close to $0$ as we want.
But by Poincar\'e's inequality $|m_{i,1}-m_{i,2}| \leq C \rho^{1-n} \mu_{W_i}(B(x,\rho))$, so
$\H^{n-1}(\d^\ast \cap B(x,\rho)) \geq  \mu_{W_i}(B(x,\rho)) \geq C^{-1} \rho^{n-1}$,
 in contradiction with
\eqref{4a33}. 

We may now summarize: we took $B(x,r)$ as in the assumptions of the theorem, assumed that \eqref{4a23}
holds, and got a contradiction if $\varepsilon$ and $r_{00}$ are chosen small enough. 
So \eqref{4a23} fails, and \eqref{4a8} holds for $r \leq r_{00}$. We can easily extend this to $r_{00} < r \leq 1$, say, 
simply by using the result for $r_{00}$ and multiplying the constant $C_{ar}$ by $r_{00}^{1-n}$ (but the true invariant 
statement is rather to assume that $r \leq r_{00}$ and get the first constant $\varepsilon$).
This completes the proof of Theorem \ref{t4a1}.
\qed

Theorem \ref{t4a1} is not always precise enough, because we also wish to know that every $W_i$
is individually nice. The following statement gives some initial information about this.

\begin{thm}\label{t4a4} 
Let $\Omega$, $a$, $G$, and $J_w$ be as in Theorem \ref{t4a1}. There exists 
a constant $C_N$, which depends only on $n$, $\delta$, $\alpha$, $C_\alpha$, but also the
number of components $W_i$, such that if $\Omega' \subset \Omega$ is open, if
\begin{equation} \label{4a34}
\delta \leq a(x) \leq 1 \ \text{ for } x \in \Omega',
\end{equation}
and if $\W \in \cF_w(\Omega)$ minimizes $J_w$ in the class $\cF_w(\Omega)$, then for $1 \leq i \leq N$,
$x \in \Omega' \cap \d^\ast W_i$ and $0 < \min(1,\dist(x,\d \Omega'))$,
\begin{equation}\label{4a35}
C_{N}^{-1} r^{n-1} \leq \H^{n-1}(\d^\ast W_i \cap B(x,r)).
\end{equation}
\end{thm}

So none of the $\d^\ast_i$ can be too small locally.
We did not need to write the upper bound, because $\d^\ast$ is larger. 
\ms
\begin{rem} \label{r4a5}
Unfortunately our estimate depends on $N$. The author's bet is that the same result also holds with
a constant $C_N$ that does not depend on $N$, but were not able to prove this. Similarly, there should be an 
automatic regulation of the number of pieces, in the sense that if $\Omega$ has a reasonable shape
(so that a Poincar\'e inequality holds in $\Omega$), and $\delta \leq a(x) \leq 1$ in the whole $\Omega$, 
and $\W \in \cF_w(\Omega)$ minimizes $J_w$ in the class $\cF_w(\Omega)$, the number of pieces $W_i$ such that
$|W_i| > 0$ should be bounded by a constant $N_0$ that depends only on $n$, $\delta$, 
$\alpha$, $C_\alpha$, and geometric constants related to $\Omega$.
A proof of this would need to be more clever than the proof below, and bound the number of components that 
really interact with a given one that we want to erase.
\end{rem}

\begin{rem} \label{r4a6}
In \cite{Le}, G. Leonardi proves an analogue of Theorem \ref{t4a1} in the more general situation where
$F_w$ looks like
\begin{equation}\label{4a36}
F_w(\W) = a_{i,j} \sum_{i, j \in I, i\neq j} \H^{n-1}(\d^\ast_i \cap \d^\ast_j),
\end{equation}
where the coefficients $a_{i,j}$ are symmetric ($a_{i,j} = a_{j,i}$) and satisfy a condition 
inspired of the triangle inequality that prevents thin layers of some $W_k$ to be 
incorporated between $W_i$ and $W_j$ so that the energy of the interfaces diminishes.
His proof is like the proof of Theorem \ref{t4a1}, but with a beautiful extra argument of graph theory
whose purpose is to take into account the fact that when we pour $W_i$ into $W_\ell$, some of the old
boundary $\d^\ast_i$ becomes a part of $\d^\ast_\ell$, so that in the functional $F_w(\W)$ the corresponding 
multiplying coefficients $a_{i,j}$ are replaced by $a_{\ell,j}$. The same proof should also work with
\begin{equation}\label{4a37}
F_w(\W) =  \sum_{i, j \in I, i\neq j} \int_{\d^\ast_i \cap \d^\ast_j} a_{i,j} d\H^{n-1}(x),
\end{equation}
provided that the variable coefficients $a_{i,j}(x)$ are also H\" older continous. 
We decided not to pursue this. As far as the authors know, it is not known whether the natural variant of 
Theorem \ref{t4a4} in this context of different coefficients holds, even with constant coefficients $a_{i,j}$
and only $N=4$.
\end{rem}

The proof of Theorem \ref{t4a1} is easy (but the result is probably too weak).
We proceed as above, suppose that $x \in \d^\ast_i$ and $r_{00}$ are such that 
\begin{equation}\label{4a38}
\H^{n-1}(\d^\ast_i \cap B(x,r_0) \leq \varepsilon r_0^{n-1},
\end{equation}
just as in \eqref{4a23} but now we will allow $\varepsilon$ and $r_{00}$ to depend on $N$. 
We choose $\rho \in (\frac{9r_0}{10}, r_0)$ as above, so that \eqref{4a12} holds for $\rho$, and 
now want to apply the same argument as above, but where we just poor some part of $W_i$ alone
into some other piece $W_\ell$. We consider $D B(x, \rho)$, add a second copy $D'$ of $D$ to
$D$ as above, to make a topological sphere $\wh D$, take $W = D \cap W_i$ and glue a second
copy $W'$ to make a set of finite perimeter $\wh W \subset \wh D$, and apply Poincar\'e's
inequality in $\wh D$. We get, just as for \eqref{4a25}, that 
\begin{multline} \label{4a39}
\min(|W|, |D\sm W|) = \frac12 \min(|\wh W|, |\wh B\sm \wh W|) 
\\
\leq C \Per(\wh W; \wh B)^{\frac{n}{n-1}} = C \H^{n-1}(\d^\ast_i \wh D)^{\frac{n}{n-1}}
\leq C \varepsilon^{\frac{n}{n-1}} r_0^n.
\end{multline}
 We start a discussion with the most likely case when $|W| = |W_i \cap D| \leq \frac12 |D|$, which
actually implies that $|D|^{-1} |W_i \cap D|$ is very small. We choose a new radius $r$, 
with $r_1 < r < \rho$, such that \eqref{4a12} 
holds for $r$ and 
\begin{multline} \label{4a40}
\int_{W_i \cap \d B(x,r)} d\H^{n-1} \leq (\rho-r_1)^{-1} \int_{r_1 < r < \rho}
\int_{W_i \cap \d B(x,r)} d\H^{n-1} dr
\\
= C (\rho-r_1)^{-1} \big|(W_i \cap B(x,\rho) \sm B(x,r_1) \Big|
\leq  C r_0^{-1} |W_i \cap D|  \leq C \varepsilon^{\frac{n}{n-1}} r_0^{n-1}
\end{multline}
with the same justification as for \eqref{4a29}. Recall from Lemma \ref{l4a1} 
that modulo a $\H^{n-1}$-negligible set,  $\d^\ast_i$ is the union of the $\d^\ast_i \cap \d^\ast_j$,
$j \neq i$, so we can choose $j \neq i$ such that 
\begin{equation}\label{4a41}
\H^{n-1}(\d^\ast_i \cap \d^\ast_j  \cap B(x,r)) \geq N^{-1} \H^{n-1}(\d^\ast_i \cap B(x,r)),
\end{equation}
and we decide to pour $W_i \cap B(x,r)$ into $W_j$. In the previous construction this amounts
to taking $I_0 = \{ i \}$ and $\ell(i) = j$. Let $X$ denote the new competitor that we get; this time
we do not add anyone to $\d^\ast$ in $B = B(x,r)$, but we remove the whole set $\d^\ast_i \cap \d^\ast_j  \cap B(x,r)$,
so we save at least the right-hand side of \eqref{4a41} there. As before, we add the new piece $W_i \cap \d B$
to $\d^\ast_i$ and $\d^\ast_i$, and we do not change $\d^\ast$ on $\Omega \sm \ol B$. Altogether,
\eqref{4a14} and \eqref{4a15} now yield
\begin{equation}\label{4a42}
\H^{n-1}(\d^\ast(\W) \sm \d^\ast(\X)) 
\geq \H^{n-1}(\d^\ast_i \cap \d^\ast_j  \cap B) \geq N^{-1} \H^{n-1}(\d^\ast_i \cap B)
\end{equation}
and
\begin{multline} \label{4a43}
\H^{n-1}(\d^\ast(\X) \sm \d^\ast(\W))
\leq \H^{n-1}(\d^\ast(\X) \cap \d B)
\\
\leq \H^{n-1}(\d B \sm W_1) \leq C (a^k \varepsilon)^{\frac{n}{n-1}} r_k^{n-1}
\end{multline}
by \eqref{4a14} and \eqref{4a15}. Hence by \eqref{4a22}
\begin{multline} \label{4a44}
\H^{n-1}(\d^\ast_i \cap B(x,r_1)) 
\leq \H^{n-1}(\d^\ast_i \cap B) \leq N \H^{n-1}(\d^\ast(\W) \sm \d^\ast(\X)) 
\\
\leq C N\H^{n-1}(\d^\ast(\X) \sm \d^\ast(\W)) + C'_\gamma N r_0^\gamma r_0^{n-1}
\leq C N \varepsilon^{\frac{n}{n-1}} r_0^{n-1} + C'_\gamma N r_0^\gamma r_0^{n-1}.
\end{multline}
That is, $r_1$ satisfies \eqref{4a38} with $\varepsilon$ replaced by
$\varepsilon_1 = C \varepsilon^{\frac{n}{n-1}} N + C'_\gamma r_0^\gamma N$.
This should of course be compared to \eqref{4a31}; we see that we only had to multiply the
right-hand side by $N$.

Before we can conclude as above, we still need to say what happens when 
$|W| = |W_i \cap D| > \frac12 |D|$, i.e. when $W_i$ is the major part of $D$. 
In this case we proceed exactly as in Theorem \ref{t4a1}, i.e., choose
$r$ such that \eqref{4a29} holds (with $W_1$ replaced by $W_i$), pour all the other parts
$W_j \cap B$, $j \neq i$, into $W_i$, and get \eqref{4a31}, which is even better than \eqref{4a44}
by a factor $N$.

So we managed to prove \eqref{4a44}, and now we can copy the induction argument and the rest of the
proof of Theorem \ref{t4a1}; Theorem \ref{t4a4} follows.
\qed

\section{Clean representatives of weak minimizers}\label{S5} 

In this section we take a minimizer $\W \in \cF_w(\Omega)$ for the weak functional $J_w$, for instance,
and prove that we can clean it and get an equivalent $N$-tuple $\W^\sharp$, where the $W_i^\sharp$ 
are open. This will be used later to prove that the weak and strong functionals are equivalent, and the argument
is rather standard in the business of weak functionals.
We keep the same notation as above concerning
the reduced boundaries $\d^\ast_i = \d^\ast W_i$ and $\d = \cup_i \d^\ast_i \subset \Omega$.
Most of the discussion in this section does not use the form of the functional, but only 
local Ahlfors regular estimates for $\d^\ast$ like \eqref{4a8}.

\begin{pro}\label{t5a1} 
Let $\W \in \cF_w(\Omega)$ be such that for every compact set $K \subset \Omega$, 
there is a constant $C_K \geq 1$ such that for $i\in I = \{ 1,2,\ldots , N \}$,
\begin{equation}\label{5a1}
C_K^{-1} r^{n-1} \leq \H^{n-1}(\d^\ast_i \cap B(x,r)) \leq C_K r^{n-1}
\ \text{ for $x\in \d^\ast \cap K$ and $0 < r \leq 1$.} 
\end{equation}
Then 
\begin{equation}\label{5a2}
\H^{n-1}(\Omega \cap\ol{\d^\ast_i} \sm \d^\ast_i) = 0 
\ \text{ for } i\in I,
\end{equation}
and we can find open sets $W_i^\sharp$, $1 \leq i \leq N$, such that for $i\in I$,
\begin{equation}\label{5a3}
|W_i \sm W_i^\sharp| = |W_i^\sharp \sm W_i| = 0
\end{equation}
and
\begin{equation}\label{5a4}
\Omega \cap \d W_i^\sharp = \Omega \cap \ol{W_i^\sharp} \sm W_i^\sharp = \ol{\d^\ast_i}.
\end{equation}
\end{pro}

\ms
If we only know \eqref{5a1} globally for $\d^\ast$, we get still \eqref{5a3} and \eqref{5a4}, and
instead of \eqref{5a2} we at least have 
\begin{equation}\label{5a5}
\H^{n-1}(\Omega \cap\ol{\d^\ast} \sm \d^\ast) = 0.
\end{equation}

Notice that the conclusion holds for minimizers of the weak functional, or  even when $\W$ lies in
one of the classes of almost- and quasi-minimizers alluded to in Definition \ref{d4a3} or near \eqref{4a22}, 
since we then have the conclusions of Theorem \ref{t4a1} or \ref{t4a4}. Set
\begin{equation}\label{5a6}
\d_i = \Omega \cap \ol{\d^\ast_i}  \ \text{ and } \  \d = \Omega \cap \ol{\d^\ast} = \cup_{i \in I} \d_i;
\end{equation}
the main part of the proposition is \eqref{5a2} (or \eqref{5a5} if we only have \eqref{5a1} for $\d^\ast$),
which we check now.

We prove \eqref{5a5}; \eqref{5a2} would be the same (add an index $i$). 
It is enough to check that $\H^{n-1}(B \cap \d \sm \d^\ast) = 0$ when
$B$ is a closed ball such that $3B \subset \Omega$. Notice that then we can apply Theorem \ref{t4a1}
on $B$. Let $\varepsilon > 0$, and set $\mu = \H^{n-1}_{\vert 2B \cap \d^\ast}$. Since by \eqref{4a8} $\mu$ is a 
Radon measure, we can find a compact set $K \subset \d^\ast$ such that $\mu(2B \sm K) \leq \varepsilon$.
Then let $\eta > 0$ be small, and for each $x\in \d \cap B \sm K$, choose a small ball $B(x) = B(x, r(x))$ 
such that $0 < r(x) < \eta/2$ and $B(x) \cap K = \emptyset$. Now use the standard $5$-covering lemma of Vitali: 
there is an at most countable set $X \subset \d \cap B \sm K$ such that the balls $B(x, r(x)/5)$, $x \in X$, 
are disjoint, but the $B(x)$ cover $\d \cap B \sm K$. 
Since $\d = \Omega \cap \ol{\d^\ast}$, each $B(x, r(x)/10$ contains a point $y \in \d ^\ast$, and by \eqref{5a1} 
for $\d^\ast$,
\begin{equation}\label{5a7}
\mu(B(x, r(x)/5)) \geq \mu(B(y, r(x)/10)) = \H^{n-1}(\d^\ast \cap B(y, r(x)/10))\geq C^{-1} r(x)^{n-1}.
\end{equation}
Now
\begin{equation}\label{5a8}
\sum_{x\in X} r(x)^{n-1} \leq C \sum_{x\in X} \mu(B(x, r(x)/5)) \leq C \mu(2B\sm K) \leq C\varepsilon
\end{equation}
because the balls $B(x, r(x)/5)$ are disjoint and contained in $2B\sm K$. In terms of Hausdorff measure, 
this shows that $\H^{n-1}_{\eta}(\d \cap B \sm K) \leq C \varepsilon$ (for the definition of $\H^{n-1}_{\eta}$, see chapter 4 of \cite{Ma}).
 Since $\eta$ is as small as we want,
we get that $\H^{n-1}(\d \cap B \sm K) \leq C \varepsilon$, and since $\varepsilon$ is arbitrary, 
$\H^{n-1}(\d \cap B \sm K) = 0$; \eqref{5a5} follows.

Now define the $W_i^\sharp$, $i \in I$, by
\begin{equation}\label{5a9}
W_i^\sharp = \big\{ x\in \Omega \sm \d \, ; \, |B(x,r) \sm W_i| = 0 \text{ for some } r > 0 \big\}.
\end{equation}
These are clearly open sets, and we claim that
\begin{equation}\label{5a10}
\Omega \sm \d \text{ is the disjoint union of the $W_i^\sharp$, $i \in I$.}
\end{equation}
Indeed, let $x \in \Omega \sm \d$ be given; by definition we can find $r > 0$ such that $B = B(x,r)$ is 
contained in $\Omega$ and does not meet $\d^\ast$. Since the $W_i$ cover $\Omega$, we can find $i$ such that 
$|W_i \cap B| > 0$. Now if $|B \sm W_i| > 0$,  we get for instance that 
$\int_{x\in B}\int_{y\in B} |\1_{W_i(x)}-\1_{W_i(y)}| dxdy \geq |W_i \cap B||B \sm W_i| > 0$, and by
Poincar\'e's inequality $\mu_{W_i}(B) = \int_B |D\1_{W_i}| > 0$. But $\mu_{W_i}(B) = \H^{n-1}(\d^\ast_i \cap B)$ by definition of $\mu_{W_i}(B)$ and \eqref{4a1}, and this contradicts the fact that $B$ does not meet 
$\d^\ast_i \subset \d^\ast$. So  $x\in W_i^\sharp $, and \eqref{5a10} follows.

Next we check the equivalence \eqref{5a3}. If $x \in W_i$, then almost always (by \eqref{5a5})
$x \in \Omega \sm \d$, hence (by \eqref{5a9}) it lies in some $W_j^\sharp$. But $x$ is almost always a Lebesgue
density point for $W_i$, so $j=i$ and $x\in W_i^\sharp$. Conversely, if $x\in W_i^\sharp$, 
some whole ball $B(x,r)$ is (almost) contained in $W_i$; but almost always $x$ is a point of density of all the 
$W_j$ that contain it, and this forces $j=i$. Hence $x\in W_i$ (because $\Omega = \cup_j W_j$).

We still need to check \eqref{5a4}. Let $x \in \Omega \cap \d W_i^\sharp = \ol{W_i^\sharp} \sm W_i^\sharp$.
Every small ball $B(x,r)$ meets $W_i^\sharp$, so $|W_i \cap B(x,r)| > 0$ (because $W_i^\sharp$ is open).
On the other hand, $x \notin W_i^\sharp$ hence $|B(x,r) \sm W_i| > 0$ so, by the same proof using 
the Poincar\'e inequality as before, $B(x,r)$ meets $\d^\ast_i$. So $\d W_i^\sharp \subset \d_i$.
Conversely, if $x \in \d_i$, every ball $B(x,r)$ meets $\d^\ast_i$. If $y$ is any point of $B(x,r)\cap\d^\ast_i$
and $D$ is a tiny ball centered at $y$, then $|W_i \cap D| > 0$ and $|D\cap W_i| > 0$ by definition
of $\d^\ast_i$ (see \eqref{4a3}). Then also $|W_i^\sharp \cap D| > 0$ and $|D\cap W_i^\sharp| > 0$ by
\eqref{5a3}, and $x \in \d W_i^\sharp$. Proposition \ref{t5a1} follows.
\qed

\ms
The $N$-tuple $\W^\sharp$ does not exactly lie is our set $\cF$ of strong partitions, because
it is not a partition. But we can easily find a true partition $\wh\W \in \cF$, with the property that
for $i\in I$,
\begin{equation}\label{5a11}
W_i^\sharp \subset \wh W_i \subset \Omega \cap\ol{W_i^\sharp}
\end{equation}
and then, because $\Omega \cap \d W_i^\sharp \subset \d$ has no interior and by \eqref{5a4}
\begin{equation}\label{5a12}
\Omega \cap \d \wh W_i = \Omega \cap \d W_i^\sharp = \Omega \cap \ol{\d^\ast_i}.
\end{equation}
Now we can compute the strong functional. Of course $G(\wh\W)$ is the same for $J_w$ and $J_s$, 
so we are left with
\begin{multline}\label{5a13}
F_s(\wh\W) = \sum_{i=1}^N \int_{\Omega \cap \d \wh W_i} a(x) \d\H^{n-1}(x)
= \sum_{i=1}^N \int_{\Omega \cap\ol{\d^\ast_i}} a(x) \d\H^{n-1}(x)
\\
= \sum_{i=1}^N \int_{\Omega \cap \d^\ast_i} a(x) \d\H^{n-1}(x) 
= F_w(\wh\W),
\end{multline}
by \eqref{3a3}, \eqref{5a12}, and \eqref{5a2}. The authors are not sure that there may be circumstances
where one can obtain \eqref{5a1} globally for $\d^\ast$ but not separately for each $\d^\ast_i$, but if this is the
case we can also use a variant of $F_s$, namely
\begin{equation} \label{5a14}
F_s^\sharp(\W) = 2\int_{\d(\W)} a(x) \d\H^{n-1}(x),
\end{equation}
where $\d(\W) = \Omega \cap (\cup_i \d W_i)$ as above. For this one, even if we only have \eqref{5a1} for $\d^\ast$,
\begin{equation} \label{5a15}
F_s^\sharp(\wh\W) = 2\int_{\d^\ast} a(x) \d\H^{n-1}(x) =  F_w(\wh\W),
\end{equation}
by \eqref{5a5} and \eqref{4a5}.

We are now ready to say that the weak and strong functionals are equivalent.

\begin{cor}\label{c5a2} 
Let $\Omega$ be a bounded domain, $a$ a continuous positive function on $\Omega$, and $G$ satisfy the 
H\"older condition \eqref{4a6}. Then the functional $J_w$ of \eqref{3a7} has a
minimizer on $\cF_w$, and the associated cleaner partition $\wh\W$ constructed above minimizes 
the functionals $J_s$ and $J_s^\sharp$ on $\cF$. In particular,
\begin{equation} \label{5a16}
\inf_{\W \in \cF} J_s(\W) = \inf_{\W \in \cF} J_s^\sharp(\W) =  \inf_{\W \in \cF_w} J_w(\W) 
= J_w(\wh\W).
\end{equation}
Also, if $\W_0$ is a minimizer of $J_s(\W)$ or $J_s^\sharp(\W)$ in $\cF$,
then it is also a minimizer of $J_w$ in $\cF_w$.
\end{cor}

Notice first that the assumptions of Theorem \ref{t3a1} are satisfied, so $J_w$ has a minimizer
in $\cF_w$. In addition, the assumptions of Theorem \ref{t4a4} are also satisfied, except for the fact that
$a(x) \leq 1$ on $\Omega$, which is not needed for the local Ahlfors estimate (see the remarks below the statement
of Theorem~\ref{t4a1}, or simply follow the proof), so the $\d^\ast_i$ satisfy the assumptions of
Proposition~\ref{t5a1}, and we can construct  the clean representatives
$\W^\sharp$ and $\wh\W$ with the properties \eqref{5a14} and \eqref{5a15}.
Then $J_s^\sharp(\wh\W) = J_s(\wh\W) = J_w^\sharp(\wh\W)$.

Now it should be observed that for any $\W \in \cP$ such that $J_s(\W)$ or $J_s^\sharp(\W)$ is finite,
the $W_i$ have locally finite perimeters, and $J_w(\W) \leq J_s(\W), J_s^\sharp(\W)$. This is because
for any set $A$ such that $\H^{n-1}(A) < +\infty$, $A$ has a finite perimeter and $\d^\ast A \subset \d A$;
then we can use the definitions and \eqref{4a5}. Hence 
$\inf_{\W \in \cF} J_s(\W)$ and $\inf_{\W \in \cF} J_s^\sharp(\W)$ are at least as large
as $\inf_{\W \in \cF_w} J_w(\W) = J_w(\wh\W)$. The rest of the corollary follows.
\qed

\section{Uniform rectifiability and isoperimetry}
\label{S6} 

In this section we consider a minimizer $\W$ of the functional $F_w$ on $\Omega$,
and assume as in Theorems \ref{t4a1} and \ref{t4a4} that $G$ is H\"older continuous, as in \eqref{4a6},
and that on the open set $\Omega' \subset \Omega$, we have the bounds \eqref{4a7} on $a$.

We continue the description of the free boundary.
 In fact, we shall find it more convenient to use the 
equivalent open sets 
$\W_i^\sharp$ of Proposition \ref{t5a1} (defined by \eqref{5a9}), and the closed boundaries
\begin{equation}\label{6a1}
\d_i = \Omega \cap \ol{\d^\ast_i} = \Omega \cap \d W_i^\sharp 
 \ \text{ and } \  \d = \Omega \cap \ol{\d^\ast} = \cup_{i \in I} \d_i;
\end{equation}
of \eqref{5a6} (see \eqref{5a3}). Recall from \eqref{5a2} that 
$\H^{n-1}(\Omega \cap\ol{\d^\ast_i} \sm \d^\ast_i) = 0$,  so we are not adding any mass when we replace 
$\d^\ast_i$ with $\d_i$, but having a closed set is nicer. Because of Corollary \ref{c5a2}, we could also have 
considered a minimizer of $J_s$ or $J_s^\sharp$ in $\cF$, and cleaned it as we have cleaned weak minimizers
in Section~\ref{S5}.

Recall that Theorems \ref{t4a1} says that $\d$ is a locally Ahlfors regular set in $\Omega'$, and 
Theorems \ref{t4a4} says that the individual $\d_i$ are also locally Ahlfors regular sets in $\Omega'$,
but with bounds that depend on the number $N$ of pieces too. Our next step is to prove that these sets
satisfy S. Semmes' Condition B, and hence are locally uniformly rectifiable. 

\begin{thm}\label{t6a1} 
Let $\Omega \subset \R^n$,  $\W \in \cF_w(\Omega)$, $\d$, and $\Omega'$ be as above.
There exists a constant $C_1 = C_1(n,\alpha, C_\alpha, \delta) \geq 1$ such that
for $x \in \d \cap \Omega'$ and $0 < r \leq \min(1, \dist(x,\d \Omega'))$, we can find two indices $i, j \in I$,
with $i \neq j$, and points $x_i \in W_i^\sharp$ and $x_j \in W_j^\sharp$, such that 
\begin{multline} \label{6a2}
B_i = B(x_i, C_1^{-1} r) \subset W_i^\sharp \cap B(x,r)
\\ \text{ and } \ 
B_j = B(x_j, C_1^{-1} r) \subset W_j^\sharp \cap B(x,r).
\end{multline}
Moreover, there exists a constant $C_2 = C_2(n,\alpha, C_\alpha, \delta,N) \geq 1$ such that
for $x \in \d_i \cap \Omega'$, with $i \in I$, and $0 < r \leq \min(1, \dist(x,\d \Omega'))$, 
we can find $x_i \in W_i^\sharp$  such that 
\begin{equation} \label{6a3}
B(x_i, C_2^{-1} r) \subset W_i^\sharp \cap B(x,r)
\end{equation}
\end{thm}

As before, it could be that \eqref{6a3} also holds with a constant $C_2$ that does not depend on $N$, 
but we shall not prove this. Condition B was introduced in \cite{Se}, and by now we have various proofs 
of the fact that for locally Ahlfors regular sets, it implies the local uniform rectifiability of $\d$ 
(or of the individual $\d_i$, with a worse constant), and even that $\d$ or the $\d_i$ contain 
``big pieces of Lipschitz graphs''. 
See \cite{DS} for all the information that the reader could want on uniform rectifiability; here we shall only 
consider this as good news, and shall not comment further.

The proof will follow the same line as in \cite{DS2}. We start with the result with the two balls. 
Notice that we need $N \geq 2$ for \eqref{6a4} to hold, but this is all right because 
if $N=1$, $\Omega = W_1$ and $\d$ is empty. 

Let $(x,r)$ be as in the statement; without loss of generality, we can assume that 
$|B(x,r/2) \cap W_1| \geq |B(x,r/2) \cap W_2| \geq |B(x,r/2) \cap W_i|$ for $i > 2$. 
Let us first show that 
\begin{equation}\label{6a4}
|B(x,r/2) \sm W_1| \geq \varepsilon r^n
\end{equation}
for some $\varepsilon > 0$ that depends on $n$, $\alpha$, $C_\alpha$, and $\delta$. 
If not, pick $\rho \in (r/4, r/2)$, as we chose $r$ near \eqref{4a29}, such that \eqref{4a12} holds for 
$\rho$ and 
\begin{equation}\label{6a5}
\H^{n-1}(\d B(x,\rho) \sm W_1) \leq C r^{-1} |B(x,r/2) \sm W_1| \leq C \varepsilon r^{n-1}.
\end{equation}
Then try the competitor $\X$ such that $X_1 = W_1 \cup B(x,\rho)$ and
$X_i = W_i \sm B(x,\rho)$ for $i > 1$. That is, pour the $W_i \cap B(x,\rho)$ into $W_1$.
The accounting is the same as near \eqref{4a30}: we save $\d^\ast \cap B(x,\rho)$, 
but we may add $B(x,\rho) \sm W_1$ to $\d^\ast$, and also add 
\begin{equation}\label{6b5}
|G(\W) - G(\X)|\leq C_\alpha \dist(\W,\X)^\alpha = C_\alpha 2^\alpha |W_i \cap B(x,\rho)|^\alpha
\leq C C_\alpha \varepsilon^\alpha r^{n\alpha}
\end{equation}
by \eqref{4a8}. The minimality of $\W$ yields
\begin{multline}\label{6a6}
\H^{n-1}(\d^\ast \cap B(x,\rho))) \leq C \delta^{-1} \H^{n-1}(\d B(x,\rho) \sm W_1) 
+ C C_\alpha \delta^{-1} \varepsilon^\alpha r^{n\alpha}
\\\leq C \delta^{-1} \varepsilon^\alpha r^{n-1}
\end{multline}
because $\alpha > \frac{n-1}{n}$.
If $\varepsilon$ is small enough, this contradicts Theorem \ref{t4a1}, so \eqref{6a4} holds.

We also need to know that $W_1$ is not too small, i.e., that 
\begin{equation}\label{6a7}
|B(x,r/2) \cap W_1^\sharp| = |B(x,r/2) \cap W_1| \geq \varepsilon r^n,
\end{equation}
where the first part comes from \eqref{5a3}.
If not, then for each $i$, $|B(x,r/2) \cap W_i| \leq \varepsilon r^n$, 
and by Poincar\'e's inequality on the double of $B(x,r/2)$, and as \eqref{4a26}, we have that for 
for each $i \in I$,
\begin{multline} \label{6a8}
|B(x,r/2) \cap W_i| = \min(|B(x,r/2) \cap W_i|, B(x,r/2) \sm W_i|) 
\\ \leq C \Per(W_i,B(x,r/2)) ^{\frac{n}{n-1}}
\leq C  \H^{n-1}(\d_i^\ast \cap B(x,r/2))^{\frac{n}{n-1}}
\end{multline}
and hence 
\begin{multline} \label{6a9}
\H^{n-1}(\d_i^\ast \cap B(x,r/2)) \geq C^{-1} |B(x,r/2) \cap W_i|^{\frac{n-1}{n}}
\\= C^{-1} |B(x,r/2) \cap W_i| \, |B(x,r/2) \cap W_i|^{-\frac{1}{n}}
\geq C^{-1} |B(x,r/2) \cap W_i| \varepsilon^{-\frac{1}{n}} r^{-1}.
\end{multline}
We sum this over $i \in I$, and get that
\begin{equation} \label{6a10}
\sum_{i \in I} \H^{n-1}(\d_i^\ast \cap B(x,r/2)) \geq C^{-1} \varepsilon^{-\frac{1}{n}} r^{-1} |B(x,r/2)| 
\geq C^{-1} \varepsilon^{-\frac{1}{n}} r^{n-1}.
\end{equation}
By Lemma \ref{l4a1}, almost every point of $\d^\ast$ lies in exactly two sets $\d^\ast_i$, so 
\begin{equation} \label{6a11}
\H^{n-1}(\d^\ast \cap B(x,r/2)) \geq C^{-1} \varepsilon^{-\frac{1}{n}} r^{n-1},
\end{equation}
which contradicts the upper bound in \eqref{4a8} if $\varepsilon$ is small enough. So \eqref{6a7}
holds too.

Let $\tau >0$ be very small (even smaller than $\varepsilon$), to be chosen soon, and set
\begin{equation} \label{6a12}
Z = \big\{ y \in B(x,r/2) \, ; \, \dist(y, \d) \leq \tau r \big\};
\end{equation}
we want to show that 
\begin{equation} \label{6a13}
|Z| \leq C \tau r^{n},
\end{equation}
where $C$ depends on $n$, $\alpha$, $C_\alpha$, and $\delta$, and for this we choose a maximal family
$\{ z_k \}$,  $k \in S$,
 of points of $Z$ at mutual distances larger than $4\tau r$. Since the balls
$B(z_k, 2\tau r)$ are disjoint, and each one contains a ball of radius $\tau r$ centered at a point of 
$B(x,3r/4)$, we see that the number $M = \# S$  of balls is such that
\begin{eqnarray} \label{6a14}
M &\leq& C \sum_{k \in K} (\tau r)^{1-n} \H^{n-1}(\d^\ast \cap B(z_k,2\tau r))
\nn\\
&\leq& C (\tau r)^{1-n} \H^{n-1}(\d^\ast \cap \big[\cup_{k} B(z_k,2\tau r)\big])
\nn \\
&\leq& C (\tau r)^{1-n} \H^{n-1}(\d^\ast \cap B(x,r))
\leq C (\tau r)^{1-n} r^{n-1} = C \tau^{1-n}.
\end{eqnarray}
by \eqref{4a8}. On the other hand, by maximality of $S$, the balls $B(z_k, 5 \tau r_k)$ cover $Z$, so
$|Z| \leq C M (\tau r)^n \leq C \rho r^n$, as announced in \eqref{6a13}.

If $\tau$ is small compared to $\varepsilon$, \eqref{6a7} implies that $B(x,r/2) \cap W_1^\sharp$
is not contained in $Z$, so we can find $x_1 \in  B(x,r/2) \cap W_1$ such that $\dist(x_1, \d^\ast) > \tau r$.
Then $B_1 = B(x_1, \rho r) \subset W_1^\sharp \cap B(x,r)$ (see \eqref{5a4}).
Similarly, 
\begin{equation} \label{6a15}
\sum_{i > 1} |B(x,r/2) \cap W_i^\sharp| = \sum_{i > 1} |B(x,r/2) \cap W_i|
= |B(x,r/2) \sm W_1| \geq \varepsilon r^n
\end{equation}
by \eqref{6a4}, so we can find $x_2 \in B(x,r/2) \cap \big[\cup_{i > 1} W_i^\sharp\big]$
that does not lie in $Z$, and then $B_2 = B(x_2, \rho r)$ is contained in the set $W_i^\sharp$, $i \geq 2$,
that contains $x_2$. This completes the first part of the theorem, with $C_1 = \tau^{-1}$.

Now we prove the second part. Let $x \in \d_i \cap \Omega'$ and $r > 0$ be as in the statement. 
We first check that 
\begin{equation}\label{6a16}
|B(x,r/2) \cap W_i^\sharp| = |B(x,r/2) \cap W_i| \geq \varepsilon r^n,
\end{equation}
where now $\varepsilon > 0$ is allowed to depend on $N$ too.
Indeed, if not we can choose $\rho \in (r/4, r/2)$ as before, such that \eqref{4a12} holds for $\rho$ and 
\begin{equation}\label{6a17} 
\H^{n-1}(\d B(x,\rho) \cap W_i) \leq C r^{-1} |B(x,r/2) \cap W_i| \leq C \varepsilon r^{n-1}.
\end{equation}
Then try the competitor $\X$ obtained by pouring $W_i \cap B(x,\rho)$
into another component $W_j$ such that $\H^{n-1}(\d^\ast_i \cap \d^\ast_j \cap B(x,\rho))$
is largest, so that
\begin{eqnarray} \label{6a18}
\H^{n-1}(\d^\ast_i \cap \d^\ast_j \cap B(x,\rho)) &\geq& \frac{1}{N} 
\sum_{\ell \neq i} \H^{n-1}(\d^\ast_i \cap \d^\ast_\ell \cap B(x,\rho))
\nn\\
&=& \frac{1}{N} \H^{n-1}(\d^\ast_i \cap B(x,\rho))
\end{eqnarray}
where we used again the fact that almost every point of $\d^\ast$ lies in exactly two
sets $\d^\ast_j$. Now the accounting is that we save 
$\H^{n-1}(\d^\ast_i \cap \d^\ast_j \cap B(x,\rho))$ (see \eqref{4a11}), we pay 
for the additional boundary $\d B(x,\rho) \cap W_i^\ast$, and also 
\begin{eqnarray}\label{6b20}
|G(\W) - G(\X)| &\leq& C_\alpha \dist(\W,\X)^\alpha 
\nn\\
&=& C_\alpha 2^\alpha |W_i \cap B(x,\rho)|^\alpha
\leq C C_\alpha \varepsilon^\alpha r^{n\alpha}
\end{eqnarray}
The minimality of $\W$ yields
\begin{eqnarray}\label{6a19}
\H^{n-1}(\d^\ast_i \cap B(x,\rho))) 
&\leq& C N \delta^{-1} \H^{n-1}(\d B(x,\rho) \sm W_1) + C \delta^{-1} C_\alpha \varepsilon^\alpha r^{n\alpha}
\nn\\
&\leq& C N \delta^{-1} \varepsilon^\alpha r^{n-1}.
\end{eqnarray}
If $\varepsilon$ is small enough, this contradicts Theorem \ref{t4a4}, so \eqref{6a16} holds.
And again this implies the existence of $x_i$ such that \eqref{6a3} holds, because \eqref{6a16} says
$B(x,r/2) \cap W_i^\sharp$ cannot be contained in the set $Z$ of \eqref{6a12}, if $\tau$ is chosen small enough
(now depending on $N$ too).
Theorem \ref{t6a1} follows.
\qed

\ms
We continue the description of the minimizer $\W$ by saying that the open sets $W_i^\sharp$ have a rather 
round shape locally. We already know from Condition $B$ that they are not too small, since if $x \in \d_i \cap \Omega'$ 
lies at distance $r$ from $\d \Omega'$, we can find a ball of size $C_2^{-1} \min(1, r)$ near $x$
that is contained in $W_i$. Now we want to say that near $x$, $W_i^\sharp$ is also reasonably round,
which we will express in terms of local isoperimetric inequalities in $W_i^\sharp$.
We shall use the notation of \eqref{3a1} for local perimeters of sets.

\begin{thm}\label{t6a2} 
Let $\Omega \subset \R^n$,  $\W \in \cF_w(\Omega)$, $\d$, and $\Omega'$ be as above.
There exists constants $v_0 = v_0(n,\alpha, C_\alpha, \delta, N) > 0$ and  
$C_3 = C_1(n,\alpha, C_\alpha, \delta, N) \geq 1$ such that if $i\in I$ and
$Z \subset W_i^\sharp$ is a set of finite perimeter such that  $Z \subset\subset \Omega'$ and 
$|Z| \leq v_0$, then
\begin{equation}\label{6a20}
|Z| \leq C_3 \Per(Z; W_i^\sharp)^{\frac{n}{n-1}}.
\end{equation}
 \end{thm}

We cannot allow $|Z|$ to be too large, because $Z$ could be the whole set $W_i^\sharp$, which is allowed
to be large by the $G$-term of the functional..
Of course the usual isoperimetric inequality says that 
\begin{equation}\label{6a21}
|Z| \leq C \Per(Z)^{\frac{n}{n-1}} = C \H^{n-1}(\d^\ast Z)^{\frac{n}{n-1}},
\end{equation}
so if \eqref{6a20} fails, this implies that 
$\Per(Z; W_i^\sharp) = \H^{n-1}(\d^\ast Z \cap W_i^\sharp)$ is very small compared to
$\H^{n-1}(\d^\ast Z)$, or equivalently that most of $\d^\ast Z$ lies on $\d W_i^\sharp$
(recall that $\d^\ast Z \subset \ol{W_i^\sharp}$ because $Z \subset W_i^\sharp$).

So we assume that \eqref{6a20} fails, and construct a competitor $\X$ by pouring $Z$ into some other set 
$W_j$. As usual, we choose $j$ so that 
\begin{eqnarray}\label{6a24}
\H^{n-1}(\d^\ast_i \cap \d^\ast_j \cap \d^\ast Z \cap \d W_i^\sharp)
&\geq& \frac{1}{N} \sum_{j \neq i} \H^{n-1}(\d^\ast_i \cap  \d^\ast_j \cap \d^\ast Z \cap \d W_i^\sharp)
\nn\\
&\geq& \frac{1}{N} \H^{n-1}(\d^\ast Z \cap \d W_i^\sharp)
\end{eqnarray}
because \eqref{5a4} says that $\d^\ast Z \cap \d W_i^\sharp \subset \Omega \cap \d W_i^\sharp
\subset \Omega \cap \ol{\d_i^\ast}$, \eqref{5a2} says that almost every point of $\Omega \cap \ol{\d_i^\ast}$ lies on 
$\d^\ast_i$, and Lemma \ref{l4a1} says that almost every point of $\d^\ast_i$ lies on some other $\d^\ast_j$.

When we pour $Z$ into $W_j$, we erase $\d^\ast_i \cap \d^\ast_j \cap \d^\ast Z \cap \d W_i^\sharp$
from the reduced boundary of both $W_i$ and $W_j$ (see the proof of \eqref{4a11}).
 We also add
the set  $\d^\ast Z \cap W_i$ 
 to the reduced boundaries of $W_j$ and $W_i$, and as usual
\begin{equation}\label{6a25}
|G(\W) - G(\X)| \leq C_\alpha \dist(\W,\X)^\alpha 
= C_\alpha 2^\alpha |Z|^\alpha
\end{equation}
and the minimality of $\W$ yields
\begin{equation}\label{6a26}
\H^{n-1}(\d^\ast_i \cap \d^\ast_j \cap \d^\ast Z \cap \d W_i^\sharp)
\leq C \delta^{-1} \H^{n-1}(\d^\ast Z \cap W_i^\sharp) + C \delta^{-1} C_\alpha  |Z|^\alpha
\end{equation}
Recall that 
\begin{equation}\label{6a27}
H^{n-1}(\d^\ast Z \cap W_i^\sharp) = \Per(Z; W_i^\sharp) \leq (C_3^{-1} |Z|)^{\frac{n-1}{n}}
\end{equation}
because \eqref{6a20} fails, so 
\begin{equation}\label{6a28}
\H^{n-1}(\d^\ast_i \cap \d^\ast_j \cap \d^\ast Z \cap \d W_i^\sharp)
\leq C \delta^{-1} (C_3^{-1} |Z|)^{\frac{n-1}{n}} + C \delta^{-1} C_\alpha  |Z|^\alpha
\end{equation}
In the other direction, recall that $\d^\ast Z \cap \d W_i^\sharp$ accounts for at least half of 
$\d^\ast Z$, so altogether by \eqref{6a21} and \eqref{6a24}
\begin{equation}\label{6a29}
|Z|^{\frac{n-1}{n}} \leq C  
\H^{n-1}(\d^\ast Z) \leq CN \delta^{-1} (C_3^{-1} |Z|)^{\frac{n-1}{n}} + CN \delta^{-1} C_\alpha  |Z|^\alpha
\end{equation}
We are allowed to take $C_3$ so large, depending also on $N$ and $\delta$, that the first term of the right-hand side gets eaten by the left-hand side. Also $\alpha > \frac{n-1}{n}$, so the inequality fails if $|Z|$ is too small.
This contradiction completes the proof of Theorem \ref{t6a2}.
\qed

The conclusion of Theorem \ref{t6a2}, the fact that each $W_i^\sharp$ is a domain of local isoperimetry,
is really a regularity result. We claim that, together with the  previous
 theorems \ref{t4a1} and \ref{t6a1} (local Ahlfors regularity and Condition B)
  for $\d_i$ it implies that each $W_i^\sharp$ is a local John domain. This means that every point
of $W_i^\sharp$ that is sufficiently close to $\d_i$ (compared to the distance to the usual $\Omega'$) can be connected to a Condition B ball (like $B(x_i, C_2^{-1} r)$ in \eqref{6a3}) by a thick path. The corresponding global result
was shown in \cite{DS2}, Theorem 6.1; the context is similar, but a little different because of the specific situation, but
the proof of \cite{DS2} can most probably be adapted to the local situation here (with estimates that may deteriorate when we look at open sets $\Omega' \subset \Omega$ that  tend to $\Omega$), to justify the fact that we do not give details here.

\section{Almost minimal sets}
\label{S7} 

In this section we say that if we add to the traditional assumption that there exists $\delta > 0$,
that depends on a domain $\Omega' \subset \Omega$, such that $\delta \leq a(x) \leq 1$ on $\Omega'$,
that $a$ is also H\"older continuous on $\Omega'$, then in addition to the quasiminimality properties that we have been using so far, we get that our free boundary $\d = \Omega \cap \ol{\d^\ast} = \Omega \cap \cup_{i \in I} \d^\ast W_i$
is an almost minimal set, and enjoys more regularity properties. The point of this paper is not to study these sets in general, so we shall content ourselves with a verification of the almost minimality property, plus two words about the
additional properties that one can get.

\begin{defi}\label{d7a1} 
Let $\Omega' \subset \R^n$ be an open set and $E \subset \Omega$ be closed in $\Omega'$.
We say that $E$ is an almost minimal set in $\Omega'$, with the gauge function $h$, when for each ball
$B = B(x,r) \subset\subset \Omega'$ and each Lipschitz mapping $\varphi : E \to \Omega'$ such that
\begin{equation} \label{7a1}
\varphi(y) = y \text{ for } y \in E \sm B(y,r)
\end{equation}
and 
\begin{equation} \label{7a2}
\varphi(y) \in \ol B(x,r) \text{ for } y \in E \cap B(x,r),
\end{equation}
we have that 
 \begin{equation} \label{7a3}
\H^{n-1}(E \cap \ol B(x,r)) \leq \H^{n-1}(\varphi(E) \cap \ol B(x,r)) + h(r) r^{n-1}.
\end{equation}
\end{defi}

Usually we only require the gauge function $h$ to be nondecreasing, 
such that  $\lim_{r \to 0^+} h(r) = 0$,
 and often at least a Dini condition
at the origin; here we will content ourselves with $h$ such that $h(r) \leq C r^{\beta}$
for some $C \geq 0$ and $\beta > 0$. We tried to make the definition as simple as possible; minor variants
would exist. The notion is due to Almgren \cite{Al}, and is not restricted to sets of co-dimension $1$ as here.

We could have demanded the existence of a one-parameter family of Lipschitz mappings 
$\varphi_t$, $0 \leq t \leq 1$, that connect the identity $\varphi_0$ to $\varphi_1 = \varphi$, 
but since we decided to take $\varphi(y) = y$ outside of $B(x,r)$, and $\varphi(y) \in \ol B(x,r)$ for $y \in E \cap B(x,r)$, the intermediate mapping defined by $\varphi_t(y) = t \varphi(y) + (1-t) y$ does the job.

\begin{pro}\label{t7a2} 
Let $\Omega$, $\cF_\omega$, and $J_w = G + F_w$ be as above, and let $\W$ be a minimizer of $\cF$ in $\cF_\omega$.
Assume that $\Omega$ is bounded, that $G$ satisfies the H\"older condition \eqref{4a6} and that on the open set
$\Omega' \subset \Omega$, we both have \eqref{4a34} and that 
\begin{equation} \label{7a4}
|a(x)-a(y)| \leq C_\beta |x-y|^{\beta}
\end{equation}
for some $\beta > 0$ and some constant $C_\beta \geq 0$. Then 
$\d = \Omega \cap \ol{\d^\ast (\W)}$ is an almost minimal set in $\Omega'$, with the gauge
$h(r) = C r^{\gamma}$, where 
$\gamma = \min(\beta, \alpha n - n +1)$ and a constant that depends only on $n$, $\delta$, $\alpha$, $\beta$,
$C_\alpha$, and $C_\beta$.
\end{pro}

Indeed, let $\W$ be as in the statement, set $E = \d = \Omega \cap \ol{\d^\ast (\W)}$, and let
$B = B(x,r)$ and $\varphi$ be as in Definition \ref{d7a1}. We want to construct a competitor $\X$,
and as usual we keep $\X = \W$ on $\Omega \sm B$. We need to attribute the various pieces of 
$B \sm \varphi(E)$ to pieces $W_i$.

Set $F = \varphi(E)$, and denote by $\{ U_j \}$, $j \in J$, the collection of connected components
of $B \sm F$. Also set $A_i = W_i^\sharp \cap \d B$, $i \in I$;
this is an open subset of $\d B$ (because $W_i^\sharp$ is open), but it may be empty. 
If $A_i$ is not empty, we denote by $\{ A_{i,k} \}$, $k \in K(i)$, the collection of its 
connected components. For each pair $(i,k)$, there is a component $H_{j(i,k)}$ such that, for every $y \in A_{i,k}$, 
all the points of $B$ that lie close enough to $y$ belong to $H_{j(i,k)}$. We want to take
\begin{equation} \label{7a5}
X_i = (W_i \sm B) \cup \bigcup_{i \in I} \Big(\bigcup_{k \in K(i)} H_{j(i,k)} \Big)
\end{equation}
for $i \in I$, and the main verification that we have to make is that the $X_i$ are disjoint,
i.e., since the $H_{j(i,k)}$ are contained in $B$, that for a single component $H_j$, we cannot 
have $j = j(i,k)$ and $j = j(i',k')$ unless $i = i'$.

So we give ourselves pairs $(i,k)$ and $(i',k')$, with $i \neq i'$. Thus $A_{i,k}$ and $A_{i',k'}$
lie in different sets $W_i^\sharp$ and $W_{i'}^\sharp$, and since $\d$ is the union of the boundaries of
the $W_i^\sharp$, we get that
\begin{equation} \label{7a6}
E = \d \text{ separates $A_{i,k}$ from $A_{i',k'}$ in } \Omega.
\end{equation}
It follows easily that 
\begin{equation} \label{7a7}
E \cap \ol B = \d \text{ separates $A_{i,k}$ from $A_{i',k'}$ in } \ol B
\end{equation}
(otherwise there is a path in $\ol B \sm E$ that connects a point of $A_{i,k}$ to a point
of $A_{i',k'}$, and the same path contradicts \eqref{7a6}. 

Now extend the mapping $\varphi$ to $E \cup (\Omega \sm B)$ by taking $\varphi(y) = y$
on $\Omega \sm B$; it is easy to see that the extension is still continuous  (although one may find examples where it is not Lipschitz), and we can extend it to a continuous mapping on $\Omega$, which is the identity on $\Omega \sm B$
and maps $\ol B$ into $\ol B$. Finally, we can define a continuous deformation $\{ \varphi_t \}$ by
$\varphi_t(y) = t \varphi(y) + (1-t) y$ that interpolates. Notice that $\varphi_t(\ol B) \leq \ol B$
and $\varphi_t(y) = y$ on $\d B$. We claim that because of this, 
\begin{equation} \label{7a8}
F \cap \ol B = \varphi_1(E \cap \ol B) \text{ separates $A_{i,k}$ from $A_{i',k'}$ in } \ol B.
\end{equation}
See \cite{Du}, 7-XVII-4.3. But since our connected component $H_{j}$ of $B \sm F$ would have access
to both $A_{i,k}$ and $A_{i',k'}$, we get a contradiction.

So the sets $X_i$ are disjoint. They do not need to cover, but we can put the part of $B$ that is not covered in, say 
$W_1$. By construction, the total boundary  $\d(\X) = \Omega \cap \cup_i X_i$ is contained  
in $(\d \sm B) \cup F$, and now we can compare. As usual, since we only modify the sets in $B$,
we get that 
\begin{equation}\label{7a9}
|G(\W) - G(\X)|\leq C_\alpha \dist(\W,\X)^{\alpha}
\leq 2 C_\alpha |B(x,r)|^\alpha \leq C C_\alpha r^{n\alpha}
\end{equation}
by \eqref{4a8}. Observe that
\begin{equation} \label{7a10}
\Big|a(x) \H^{n-1}(\d^\ast(\W) \cap \ol B) - \int_{\d^\ast(\W) \cap \ol B} a(y) d\H^{n-1}(y)  \Big|
\leq C_\beta r^\beta \H^{n-1}(\d^\ast(\W) \cap \ol B)
\end{equation}
and similarly for $\X$. Hence
\begin{equation} \label{7a11}
a(x) \big[\H^{n-1}(\d^\ast(\W) \cap \ol B)-\H^{n-1}(\d^\ast(\X) \cap \ol B) \big]
\leq \Delta_1 + \Delta_2, 
\end{equation}
with 
\begin{multline} \label{7a12}
\Delta_1 = \int_{\d^\ast (\W) \cap \ol B} a(y) d\H^{n-1}(y)
- \int_{\d^\ast (\X) \cap \ol B} a(y) d\H^{n-1}(y)
\\
= F_w(\W) - F_w(\X) = J_w(\W) - J_w(\X) - G(\W) + G(\X)
\\ \leq |G(\W) - G(\X)| \leq  C C_\alpha r^{n\alpha},
\end{multline}
and 
\begin{equation} \label{7a13}
\Delta_2 = C_\beta r^\beta  
\big[\H^{n-1}(\d^\ast(\W) \cap \ol B) + \H^{n-1}(\d^\ast(\X) \cap \ol B)\big].
\end{equation}
Recall that $\H^{n-1}_{\vert \d^\ast} = \H^{n-1}_{\vert \d}$, while for $\X$
we have that $\H^{n-1}_{\vert \d^\ast(\X)} \leq \H^{n-1}_{\vert \d(\X)}$. Then
\begin{multline} \label{7a14}
\H^{n-1}(E \cap \ol B(x,r)) - \H^{n-1}(\varphi(E) \cap \ol B(x,r))
\\
= \H^{n-1}(\d(\W) \cap \ol B(x,r)) - \H^{n-1}(\d(\X) \cap \ol B(x,r))
\\
\leq \H^{n-1}(\d^\ast(\W) \cap \ol B(x,r)) - \H^{n-1}(\d^\ast(\X) \cap \ol B(x,r))
\leq a(x)^{-1} (\Delta_1+\Delta_2).
\end{multline}
If $\H^{n-1}(E \cap \ol B(x,r)) \leq \H^{n-1}(\varphi(E) \cap \ol B(x,r))$, \eqref{7a3} holds
brutally. Otherwise, 
\begin{multline} \label{7a15}
\H^{n-1}(\d^\ast(\X) \cap \ol B) \leq \H^{n-1}(\d(\X) \cap \ol B) 
= \H^{n-1}(\varphi(E) \cap \ol B(x,r)) 
\\
\leq \H^{n-1}(E \cap \ol B(x,r))
= \H^{n-1}(\d^\ast(\X) \cap \ol B) \leq C r^{n-1}
\end{multline}
by the easy part of Theorem \ref{t4a1}, so $\Delta_2 \leq C C_\beta r^\beta r^{n-1}$
and now \eqref{7a14} yields 
\begin{multline} \label{7a16}
\H^{n-1}(E \cap \ol B(x,r)) 
\leq \H^{n-1}(\varphi(E) \cap \ol B(x,r)) + a(x)^{-1} (\Delta_1+\Delta_2)
\\
\leq \H^{n-1}(\varphi(E) \cap \ol B(x,r)) + C \delta^{-1}(C_\alpha  r^{n\alpha}
+ C C_\beta r^\beta r^{n-1}).
\end{multline}
Recall that $n\alpha > n-1$ by assumption, so we get \eqref{7a3}, with
$\gamma = \min(\beta, n\alpha -n +1) > 0$. This completes the proof of Proposition \ref{t7a2}.
\qed

There are not so many additional regularity properties that we can prove with the help of 
Proposition \ref{t7a2}. We could recover the Ahlfors regularity and uniform rectifiability (with a different
but more complicated proof), but at least when $n \leq3$ we can say more, and get a good local description of
$\W$.

When $n=2$, the only singularities that $\d = \d(\W)$ can have in $\Omega$ are propeller singularities, where three 
$C^{1+\sigma}$  curves of $\d$ meet at a point $x_0 \in \d$ with $2\pi/3$ angles. Then of course these curves
bound three different sets $\W_i^\sharp$. Probably a direct proof would be much faster.

When $n=3$, we also have a good description of the $2$-dimensional almost minimal sets, which
was obtained by J. Taylor \cite{Ta}. Near each of its points, $\d$ is  $C^{1+\sigma}$-equivalent to a minimal cone, 
and there are exactly three types of minimal cones: the planes, the sets $\bY$ obtained as an union of three half planes bounded by a line and that make $2\pi/3$ angles along that line, and the $\bT$-sets that are obtained as images by
an isometry of the cone over the union of the edges of a regular tetrahedron centered at the origin. In these cases too,
$\W$ is locally composed of $2$, $3$, or $4$ sets $W_i$ bounded by the faces of $\d$.

It is amusing that in the case of \cite{DFJM}, no such behavior was allowed; instead the functional managed to
either have a smooth interface, or have a small black region between the active $W_i$.

\section{Further questions}
\label{S8} 

Probably the most obvious question is whether the bounds of Theorem \ref{t4a4}, and the ensuing bounds for
the uniform rectifiability of the pieces $\d_i = \d W_i^\sharp$, really need to depend on $N$. This is connected
to the question of self-regulation: even if we start from a large enough value of $N$, do we get that 
if $\Omega$ has a nice shape and if \eqref{4a7} holds on the whole $\Omega$, the number of nontrivial 
components $W_i$ (i.e., such that $|W_i| > 0$) is bounded by a constant that depends only on 
$n$, $\alpha$, $C_\alpha$, bounds for the diameter and regularity of $\Omega$, and $\delta$?

We did not pay attention to what happens near the boundary. Probably, if we assume that \eqref{4a7} holds on 
the whole $\Omega$ and $\Omega$ has a reasonable shape (Lipschitz should be more than enough), the results
of this paper stay true at the boundary. We decided to state the results in a way that makes a discussion of what happens
near the boundary when $a(x)$ is allowed to tend to $0$ at a given rate near $\d\Omega$, but we did not pursue this.

We do not expect much difference with the results of the present paper with functions $a(x)$ (integrands) that
also depend nicely on the direction of the tangent plane to the $\d^\ast W_i$ at $x$. 
But we wish to remind the reader that in the case of ``infiltration'' (see \cite{Le}), when $a(x)$ depends also on
the pair of components $W_i$ such that $x \in \d^\ast W_i$, the interesting question of whether the Ahlfors regularity property holds in $\R^3$ and with $4$ components seems to be still open.

\section*{Acknowledgments}
This work initiated when the second author was visiting the Universit\'e de Paris-Sud (now the Universit\'e Paris-Saclay). He would like to thank Laboratoire de math\'ematiques d'Orsay for their warm hospitality. He was partially supported by the Ministry of Sciences, Research and Technology of Iran during his visit. He is also deeply indebted to all the support that Massoud Amini had for him during his graduate studies.

\Addresses

\end{document}